\documentclass{article}
\usepackage{amssymb, latexsym,amsmath,amscd}

\newtheorem{theorem}{Theorem}[section]
\newtheorem{definition}[theorem]{Definition}
\newtheorem{lemma}[theorem]{Lemma}
\newtheorem{remark}[theorem]{Remark}
\newtheorem{proposition}[theorem]{Proposition}
\newtheorem{corollary}[theorem]{Corollary}

\newcommand{\sus}{\ensuremath{^{*}}}
\newcommand{\bfi}{\bfseries\itshape}

\newcommand{\ddto}{\ensuremath{\left.\frac{d}{dt}\right|_{t=0}}}

\newcommand{\pe}{{\mathfrak p}}
\newcommand{\fg}{{\mathfrak g}}
\newcommand{\fh}{{\mathfrak h}}
\newcommand{\fm}{{\mathfrak m}}

\newcommand{\fp}{{\mathfrak q}}
\newcommand{\fq}{{\mathfrak q}}

\newcommand{\fa}{{\fg_{m_e}}}
\newcommand{\fb}{{\fm}}

\newcommand{\Ad}[1]{\mbox{Ad}_{#1}}
\newcommand{\Ads}[1]{\mbox{Ad}^*_{#1}}
\newcommand{\ads}[1]{\mbox{ad}^*_{#1}}
\newcommand{\J}{{\mathbf J}}
\newcommand{\bJ}{{\mathbf J}}
\newcommand{\cJ}{{\mathbf j}}
\newcommand{\cH}{{\mathcal H}}
\newcommand{\hxi}{{\cH - \cJ^\xi}}
\newcommand{\lp}{\left(}
\newcommand{\rp}{\right)}
\newcommand{\lcb}{\left\{}
\newcommand{\rcb}{\right\}}

\newcommand{\zero}{{\mathbf 0}}
\newcommand{\bC}{{\mathbb C}}
\newcommand{\bR}{{\mathbb R}}
\newcommand{\sym}{{G_{m_e} \cap G_\xi}}
\newcommand{\proof}{{\smallskip \noindent \bf Proof\ }}
\newcommand{\qed}{{$\blacklozenge$ \medskip}}
\newcommand{\smallfrac}[2]{{\textstyle \frac {#1} {#2}}}
\newcommand{\g}{\ensuremath{\mathfrak{g}}}

\begin{document}
\title{Bifurcation of relative equilibria in mechanical systems  
with symmetry}

\author{Pascal Chossat$^{1}$, Debra Lewis$^{2}$, Juan-Pablo  
Ortega$^{3}$, and
Tudor S. Ratiu$^{4}$}
\addtocounter{footnote}{1}
\footnotetext{Institut Nonlin\'eaire de Nice, CNRS-UNSA, 1361, route
des Lucioles, 06560 Valbonne,
France. \texttt{chossat@inln.cnrs.fr.}}
\addtocounter{footnote}{1}
\footnotetext{Department of Mathematics,
University of California, Santa Cruz, Santa Cruz, CA 95064, USA.
\texttt{lewis@math.ucsc.edu}. Research partially supported by
NSF Grant DMS-9802378.}
\addtocounter{footnote}{1}
\footnotetext{D\'epartement de Math\'ematiques,
\'Ecole Polytechnique F\'ed\'erale de Lausanne. CH--1015 Lausanne.  
Switzerland.
\texttt{Juan-Pablo.Ortega@epfl.ch}.
Research partially supported by a Fulbright/BCH
Fellowship, the US Information Agency and a Rotary
Ambassadorial Scholarship.}
\addtocounter{footnote}{1}
\footnotetext{Department of Mathematics,
University of California, Santa Cruz, Santa Cruz, CA 95064, USA, and
D\'epartement de Math\'ematiques,
\'Ecole Polytechnique F\'ed\'erale de Lausanne. CH--1015 Lausanne.  
Switzerland.
\texttt{Tudor.Ratiu@epfl.ch}. Research partially supported by
NSF Grant DMS-9802378  and FNS Grant 21-54138.98.}
\maketitle
\vspace{.5cm}
\normalsize
\begin{abstract}
The relative equilibria of a symmetric Hamiltonian dynamical system
are the critical points of the so--called \emph{augmented Hamiltonian}.
The underlying geometric structure of the system is used to  
decompose the
critical point equations and construct a collection of implicitly
defined functions and reduced equations describing the set of relative 
equilibria in a neighborhood of a given relative equilibrium. The  
structure of 
the reduced equations is studied in a few relevant situations.  
In particular, a persistence result of Lerman and  
Singer~\cite{singreleq}
is generalized to the framework of Abelian proper actions. Also, a  
Hamiltonian version of the Equivariant Branching Lemma and a study  
of bifurcations with maximal isotropy are presented. An 
elementary example is presented to illustrate the use of this approach.
\end{abstract}

\tableofcontents

\section{Introduction}

The systematic analysis of bifurcations of relative equilibria was greatly
stimulated about fifteen years ago by specific applications with 
nonconservative vector fields, namely the secondary bifurcations from 
nontrivial equilibria in hydrodynamical systems such as Couette--Taylor 
flows and Rayleigh--B\'enard  convection in a spherical shell.
The problem was attacked analytically by Chossat and Iooss
\cite{ChoIo1}, and more  qualitatively by Rand \cite{Ra}. A major
success of the analytical approach was obtained by Iooss \cite{Ioo}, who
classified the possible patterns bifurcating from a group orbit of
equilibria in a system with symmetry $O(2)$. In Moutrane \cite{Mou}, the
bifurcation of rotating waves, which are relative equilibria with a
single drift frequency, was investigated in the problem of the onset of
convection in a system with spherical symmetry.  However it was Krupa
\cite{Kru} who first developed a general theory for the  
bifurcation from
relative
equilibria. The basic tool he used was the {\bfi Invariant Slice  
Theorem} of
Palais
(see~\cite{pa, Bre}). If $G$ is a Lie group acting properly on the  
manifold
$M$, the Slice
 Theorem  establishes for each $m\in M$ an isomorphism between a
tubular neighborhood of the orbit $G\cdot m$ and the {\bfi normal  
bundle} with
base
$G\cdot m$ and fiber the normal slice $N_m$ to the tangent space to
$G\cdot m$ at $m$. It was shown by Field \cite{Fie} and then by Krupa
that any $G$--equivariant  vector field $X\in \mathfrak{X}(M)$  
admits in the
tubular neighborhood a
decomposition into the sum of two  vector fields: one,
$X_N$, defined on the normal bundle, and the other, $X_T$, defined on
the tangent bundle to $G\cdot m$. Krupa showed that the dynamical  
information,
in
particular the bifurcation properties for a parameter dependent family
of vector fields, is entirely contained in $X_N$.

The analysis of relative equilibria of conservative systems has played
a central role in the development of geometric mechanics, ranging from
the classic work of Riemann \cite{Rie} and Routh \cite{Rou2}, \cite{Rou4}
to Smale's seminal work \cite{smale}. However, the use of local singularity
theory methods, rather than explicit calculations or global topological
methods, in the analysis of conservative systems is relatively 
recent~\cite{holmes marsden, golubitsky stewart 87, ldbif, lthesis,
dellnitz melbourne marsden, marsden scheurle}.
Bifurcations of relative equilibria of Lagrangian systems and canonical 
Hamiltonian systems, i.e. Hamiltonian systems on cotangent bundles,  
with the
canonical symplectic structure and a lifted group action, have been  
studied by
Lewis \textit{et al.}~\cite{tops} and Lewis \cite{bld, qpm} using  
the reduced 
energy--momentum method developed in \cite{rem} and \cite{lbdiag}.  
This approach
uses the locked Lagrangian, the generalization of Smale's augmented  
potential
to Lagrangian systems and their Hamiltonian analogs, to  
characterize relative
equilibria as critical points of functions on the configuration  
manifold 
parameterized by elements of the algebra $\fg$ of the symmetry  
group $G$. 
A key component
of the reduced energy--momentum method is the decomposition of the  
tangent
space $T_q Q$ of the configuration manifold $Q$ at a point $q$ into the
tangent space $\fg \cdot q$ to the group orbit and an appropriate  
complement
consisting of so--called `internal' variations. The associated  
decomposition
of the relative equilibrium equations into `rigid' and `internal'  
equilibrium
conditions is analogous to the decompositions introduced by Field  \cite{Fie}
and Krupa \cite{Kru}
in the context of general equivariant vector fields. The `rigid'  
condition
can be used to determine a submanifold of `candidate relative  
equilibria'; 
imposing the remaining equilibrium conditions on this submanifold  
determines
the relative equilibria. 

Our goal is the development in the symplectic category of a  
decomposition tool
analogous to that of
Krupa that will take into account the additional structure present  
at the
kinematical level in Hamiltonian systems, without assuming all the  
ingredients
utilized in the reduced energy--momentum method. Given that many  
Hamiltonian 
systems are constructed on symplectic manifolds that are not  
cotangent bundles, 
such a tool is of much interest.
The analog  of the Invariant Slice Theorem in the symplectic  
category is given
by the
{\bfi Marle--Guillemin--Sternberg normal form}~\cite{nfm, nfgs, gs}  
(we will
refer to
it as the {\bfi MGS--normal form}) so, in principle, one could work  
as in
Krupa~\cite{Kru},
using this normal form instead of the Slice Theorem. This does not  
seem to be
the best way to proceed, since to search for relative equilibria of  
Hamiltonian 
systems one does not need to work with the Hamiltonian vector field;
there are scalar functions, the {\bfi augmented Hamiltonians},
whose critical points are precisely the relative equilibria.
Guided by Krupa's normal bundle decomposition for equivariant  
vector fields and the MGS--normal form, in
Section~\ref{Relative equilibria as critical points} we will  
construct a 
{\bfi slice mapping} with which we can decompose the critical point  
equations 
determining the relative equilibria into a system of four equations.
These split critical point equations are analyzed in
Section~\ref{the bifurcation equations} 
in a neighborhood of a given a relative equilibrium $m_e$.
Using the Implicit Function Theorem and Lyapunov--Schmidt  
reduction, we can
construct a local submanifold containing all relative equilibria  
sufficiently
near the group orbit of $m_e$. The remaining equilibrium  
conditions, called
the {\bfi reduced critical point equations}, on this submanifold can be
analyzed using standard techniques from bifurcation theory.  In  
Section~\ref{equivariance bifurcation} we study the equivariance  
properties of the reduced critical point equations.
In Section~\ref{Slice map refinements and the rigid residual  
equation} we construct
a slice mapping with respect to which one of the  reduced
critical point equations admits a simpler solution.

In Section~\ref{Persistence in Hamiltonian systems with Abelian  
symmetries}
we use the reduced critical equations and a slice mapping
constructed via the MGS--normal form to study the persistence of a  
family of
relative equilibria in a neighborhood of a non--degenerate relative  
equilibrium 
when the symmetry group of the system is Abelian. In particular, we  
generalize
to proper group actions a result from Lerman and
Singer~\cite{singreleq} originally proven for compact groups. This
result was already presented in~\cite{thesis}.

In Section~\ref{abelian g mu and equilibria} we study bifurcations  
from a 
degenerate relative equilibrium and find Hamiltonian analogs of  
bifurcation 
theorems for solutions with maximal isotropy which were first  
stated in the 
non--conservative context, namely the Equivariant Branching Lemma of 
Vanderbauwhede~\cite{vanderbauwhede 80} and Cicogna~\cite{cicogna  
81}, and 
a theorem for bifurcation of solutions with maximal isotropy group of 
complex type~\cite{melbourne 94, chossat maximal}. 

Finally we provide a simple application of the method described  
here to wave resonances.
This example is, for the most part, well known; it is included to
illustrate the implementation of the method, rather than to provide
new information.

\section{Relative equilibria as critical points}
\label{Relative equilibria as critical points}

Let $G$ be a Lie group acting smoothly on the manifold $M$ and let
$X\in\mathfrak{X}(M)$ be a smooth $G$--equivariant vector field on
$M$ with flow $F_t$. We say that the point $m_{e}\in M$ is a
{\bfi relative equilibrium\/} of the vector field $X$
if there exists an element $\xi$ of the Lie algebra $\fg$ of $G$,
called a {\bfi generator of the relative equilibrium},
such that $F_t(m_{e})=\exp(t \, \xi) \cdot m_{e}$.

Suppose now that the manifold $M$ is symplectic and
that the vector field $X$ is Hamiltonian, with associated
$G$--invariant Hamiltonian $h\in C^{\infty}(M)$. In addition,
assume that the Lie group action of $G$ on $M$ is
Hamiltonian, with associated equivariant momentum
map $\mathbf{J}:M\rightarrow\fg^{*}$. These
hypotheses imply that the
Hamiltonian vector field $X_h$ of $h$ and its flow $F_t$
are $G$--equivariant. In
this framework, the search for relative equilibria
reduces to the determination of the critical points
of a certain class of functions. Indeed, a classical result
(\cite[page 307]{fom} and~\cite[page 380]{arnold})
states that a point $m_{e}\in M$ is a relative equilibrium
of $X_h$ with generator $\xi\in\fg$ if and only
if $m_{e}$ is a critical point of the
{\bfi augmented Hamiltonian\/} $h^{\xi}:=h-\mathbf{J}^{\xi}$.
Thus, our algorithm is intended to identify 
the pairs $(m_e,\xi)\in M\times\fg$ such that 
\begin{equation}
\label{critical}
D h^{\xi}(m_e)=0.
\end{equation}
Note that if $m_e$ has nontrivial continuous symmetry, i.e. 
$\fg_{m_e} = \lcb \zeta \in \fg : \zeta_M(m_e) = 0 \rcb 
\neq \lcb 0 \rcb$, then the generator of $m_e$ is not unique.
If $\xi$ is a generator of a relative equilibrium $m_e$, then
for any $\zeta \in \fg_{m_e}$, $\xi + \zeta$ is also a 
generator.

The main goal of this section is the decomposition of the
{\bfi relative equilibrium equation}~(\ref{critical}) into a  
systems of four
equations, each defined on a space determined by the geometry of  
the problem.

Assume that $m_e$ is a relative equilibrium with generator $\xi$ and
momentum $\mu := \J(m_e)$.
Let $\fg_{m_e}$ denote the Lie algebra of the isotropy subgroup of
$m_e$ and $\fg_\mu$ the Lie algebra of the isotropy subgroup of $\mu$.
Choose complements $\fq$ of $\fg_\mu$ in $\fg$ and $\fm$ of  
$\fg_{m_e}$ in 
$\fg_\mu$, so that
\begin{equation}
\fg =\fg_\mu\oplus\fq=\fg_{m_e} \oplus \fm \oplus \fq.
\label{decompositions}
\end{equation}
The symbols $i$ and $\mathbb{P}$ with appropriate subscripts will
denote the natural injections and projections according to the
splittings~(\ref{decompositions}). For instance
$i_{\fa}:\fa\rightarrow \fg=\fa\oplus\fb\oplus\fq$
is the canonical injection of $\fa$ into $\fg$ and
$\mathbb{P}_{\fa}:\fg=\fa\oplus\fb\oplus\fq\rightarrow\fa$
extracts the $\fa$ component of any vector in $\fg$.

\begin{definition}
\label{slice mapping}
Let $V$ be a vector space and $\mathcal{U}\subset \fb^*\times V$ be  
an open
neighborhood of $(0,0)\in\fb^*\times V$.
A smooth mapping $\Psi: \mathcal{U} \subset \fb^* \times V  
\rightarrow M$
is said to be a {\bfi slice mapping} at the point $m_e\in M$
if it is a diffeomorphism onto its image satisfying the following  
conditions:
\begin{description}
\item[(SM1)] $\Psi(0,0)=m_e.$
\item[(SM2)] For any $(\eta,v)\in\mathcal{U}$ 
\begin{equation}
\label{decomposition}
T_{\Psi(\eta, v)} M = (\fm \oplus \fq) \cdot \Psi(\eta, v)
 +  T_{(\eta, v)} \Psi \cdot (\fb^* \times V).
\end{equation}
\item[(SM3)] The pullback $\cJ := \J \circ \Psi : \mathcal{U} \to  
\fg^*$ 
of the momentum map satisfies
\begin{equation}
\label{properties of F}
D \cJ(0)(\delta\eta,\delta v) 
= D \J(m_e)(T_{(0,0)}\Psi(\delta\eta, \delta v))
= \mathbb{P}_{\fb}^*\delta\eta
\end{equation}
for all $\delta \eta \in \fb^*$ and $\delta v \in V$.
\end{description}
\end{definition}

In the following proposition we explicitly construct a
slice mapping $\Psi$ at a point $m_e$ in a finite dimensional  
manifold $M$. 
\begin{proposition}
\label{build_slice}
Let $\psi:U\subset X\rightarrow M$ be a coordinate chart 
around a point $m_e$ in a finite dimensional manifold $M$ and let  
$V$ and $W$ 
be subspaces of the vector space $X$ satisfying
\begin{description}
\item[(i)] $\psi(0) = m_e$,
\item[(ii)] $T_0 \psi \cdot V$ is a complement to $\fm \cdot m_e$ in
$\ker D \J(m_e)$,
\item[(iii)] the map
\[
\begin{array}{cccc}
A:&W&\longrightarrow& (\fg_{m_e} \oplus \fq)^\circ \\
 &w&\longmapsto&D \J(m_e)T_0 \psi(w),
\end{array}
\]
is an isomorphism.
\end{description}
Let $V'$ and $W'$ be neighborhoods of the origin in $V$ and $W$  
such that
$V' \times W' \subset U$ and set 
${\mathcal U} := i_{\fm}^* (A W') \times V' \subset \fm^* \times V$.
Then the map
\[
\begin{array}{cccc}
\Psi:&\mathcal{U}\subset\fb^*\times V&\longrightarrow&M\\
&(\eta,v)&\longmapsto&\psi(v + A^{-1} \mathbb{P}_{\fb}^*\eta)
\end{array}
\]
is a slice mapping at $m_e\in M$. 
\end{proposition}
\proof
Property~\textbf{(SM1)} follows trivially from~\textbf{(i)}. 
Property~\textbf{(SM3)} follows from {\bf (ii)}, {\bf(iii)} and the 
definition of $\Psi$. 

As the first step in the proof of~\textbf{(SM2)}, we now show that
(\ref{decomposition}) holds at $(0, 0)$.
Note that {\bf (SM3)} implies that
\begin{equation}
\label{no intersection}
\mbox{ker}(D \J(m_e)) \cap T_{(0, 0)} \Psi(\fb^* \times \lcb \zero  
\rcb) 
= \lcb \zero \rcb.
\end{equation}
Combining {\bf (ii)}, (\ref{decompositions}), and (\ref{no  
intersection}), 
we obtain
\begin{eqnarray}
\dim \, T_{(0,0)}\Psi(\fb^*\times V)
&=& \dim \, \fb + \dim \, V \nonumber \\
&=& \dim \, (\ker(D \J(m_e))) \nonumber \\
&=& \dim \, P - \dim \, (\fg \cdot m_e) \nonumber \\
&=& \dim \, P - \dim \, \fb - \dim \, \fq.
\label{dim_count}
\end{eqnarray}
If $\zeta_P(m_e)=T_{(0,0)}\Psi(\delta \eta, \delta v)$,
then, since $\Psi$ satisfies {\bf (SM3)},
\[
D \J(m_e)\zeta_M(m_e)
= D \J(m_e)(T_{(0,0)}\Psi(\delta \eta, \delta v))
= \mathbb{P}_{\fb}^*\delta \eta.
\]
On the other hand, equivariance of $\J$ implies that
\[
D \J(m_e)\zeta_M(m_e) = - \ads \zeta \mu \in \fb^\circ.
\]
Hence $\ads \zeta \mu = 0$, i.e. $\zeta \in \fg_\mu$, and 
$\zeta_M(m_e) \in \fg_\mu \cdot m_e = \fm \cdot m_e$. 
Thus condition {\bf (i)} implies that $\zeta_M(m_e) = 0$.
Combining this result with (\ref{decompositions}) and (\ref{dim_count})
shows that (\ref{decomposition}) is valid at $(0, 0)$.

We now show that (\ref{decomposition}) holds for any  
$(\eta,v)\in\mathcal{U}$. 
Let $\{\xi_1,\ldots,\xi_j\}$, $\{\eta_1,\ldots,\eta_k\}$, and 
$\{v_1,\ldots,v_\ell\}$ be bases for $\fm\oplus\fp$, $\fb^*$, and $V$. 
Define the maps $u_i: {\mathcal U} \to TM$, $i = 1, \ldots, j + k +  
\ell$, by
\[
u_i(\eta, v) := \lcb \begin{array}{ll}
(\xi_i)_M(\Psi(\eta,v)) \qquad & 1 \leq i \leq j \\
T_{(\eta,v)}\Psi(\eta_{i - j}, 0) & j < i \leq j + k \\
T_{(\eta,v)}\Psi(0, v_{i - j - k}) & j + k < i \leq j + k + \ell 
\end{array} \right . .
\]
The arguments given above show that
$\lcb u_1(0, 0), \ldots, u_{j + k + \ell}(0, 0) \rcb$
is a basis for $T_{m_e}P$. Since linear independence is an open  
condition,
$\lcb u_1(\eta, v), \ldots, u_{j + k + \ell}(\eta, v) \rcb$
is a basis of $T_{\Psi(\eta,v)}M$ for $(\eta, v)$ sufficiently near  
the origin.
In particular,
\begin{eqnarray*}
T_{\Psi(\eta,v)}M
&=& \mbox{span} \lcb u_1(\eta, v), \ldots, u_{j + k + \ell}(\eta,  
v) \rcb \\
&=& {\rm span} \lcb (\xi_1)_M(\Psi(\eta,v)), \ldots,  
(\xi_j)_M(\Psi(\eta,v)) 
        \rcb \\
&& \qquad \oplus \, {\rm span} \lcb T_{(\eta,v)}\Psi(\eta_1, 0),  
\ldots, 
        T_{(\eta,v)}\Psi(\eta_k, 0) \rcb \\
&& \qquad \oplus \, {\rm span} \lcb T_{(\eta,v)}\Psi(0, v_1), \ldots, 
        T_{(\eta,v)}\Psi(0, v_\ell) \rcb  \\
&=& (\fm\oplus\fp)\cdot\Psi(\eta,v) \oplus  
T_{(\eta,v)}\Psi(\fb^*\times V),
\end{eqnarray*}
as required. 
\ \ \ $\blacklozenge$

The introduction of a slice mapping $\Psi$ allows us to decompose  
the critical point
equation~(\ref{critical})
into a system of four equations. Using property~\textbf{(SM2)}
of the slice mapping, one can conclude that the point  
$\Psi(\eta,v)\in M$ is a
relative equilibrium with generator $\xi$ if and only if
\begin{equation}
\left\{
\begin{array}{lll}
\label{split critical equation}
\text{\textbf{(RE1)}}\ \ \ i_{\fp}^* \ads \xi \cJ(\eta,v) &=& 0,\\
 & & \\
\text{\textbf{(RE2)}}\ \ \ i_{\fm}^* \ads \xi \cJ(\eta,v) &=& 0,\\
 & & \\
\text{\textbf{(RE3)}}\ \ \ D_{\fb^*} \lp \hxi \rp(\eta,v) &=& 0,\\
 & & \\
\text{\textbf{(RE4)}}\ \ \ D_V \lp \hxi \rp(\eta,v) &=& 0.
\end{array}
\right.
\end{equation}
The symbol $\mathcal{H}$ is defined by $\mathcal{H}:=h\circ\Psi$.

\begin{remark}
\normalfont
If symmetry is broken in a neighborhood of $m_e$, then
$\fg_{m_e} \cdot \Psi(\eta,v)$ is typically nontrivial. In this
case, the
first two conditions alone do not guarantee that the rigid condition
$\ads \xi j(\eta,v) = 0$ is satisfied. However, if all four conditions
are satisfied, then $D \lp \hxi \rp(\Psi(\eta, v)) = 0$; in particular,
\[
\ads \xi j(\eta,v) \cdot \eta
= - D \lp \hxi \rp(\Psi(\eta, v)) \cdot \eta_M(\Psi(\eta, v))
= 0
\]
for all $\eta \in \fg$.\ \ \ $\blacklozenge$
\end{remark}

\begin{remark}
\normalfont
Note that in order to split the critical point equation~(\ref{critical})
into~(\ref{split critical equation}),
only property~\textbf{(SM2)} of the slice mapping was utilized. As  
we shall
see in the following section, property \textbf{(SM3)} simplifies  
the analysis 
of the equations (\ref{split critical equation}).
Equations {\bf (RE1)} and {\bf (RE3)} are, by construction,  
nondegenerate
in the sense that implicit solutions to these equations always  
exist. Thus
the bifurcation analysis is carried out only on the equations obtained
by substituting the solutions of 
{\bf (RE1)} and {\bf (RE3)} into {\bf (RE2)} and {\bf (RE4)}.
$\blacklozenge$
\end{remark}

\section{The reduced critical point equations}
\label{the bifurcation equations}

In this section we start with a relative equilibrium $m_e$ with  
generator
$\xi\in\fg$
and, using the Implicit Function Theorem and 
Lyapunov--Schmidt reduction (see for instance~\cite{g1}), we derive 
a minimal set of mappings and equations determining the relative  
equilibria 
in a neighborhood of $m_e$. We will call the remaining equations the 
{\bfi reduced critical point equations}. We proceed in three steps.

\medskip

\noindent
{\bf Step 1.} Using the notation introduced in  
Definition~\ref{slice mapping},
let  $F_1: \mathcal{U} \times \fa \times \fb \times \fq \to \fp^*$  
be the
mapping
 given by
\[
F_1(\eta, v, \alpha, \beta, \gamma)
:= i_{\fp}^* \ads {\alpha + \beta + \gamma} \cJ(\eta, v),
\]
with differential
\begin{eqnarray*}
\lefteqn{D F_1(\zero) \cdot (\delta \eta, \delta v,
\delta \alpha, \delta \beta, \delta \gamma)}\\
&=& i_{\fp}^*\lp \ads {\delta \alpha + \delta \beta + \delta  
\gamma} \cJ(0, 0) 
        + \ads \xi \lp D \cJ(0, 0)(\delta\eta,\delta v) \rp \rp \\
&=& i_{\fp}^* \lp \ads {\delta \gamma} \mu+ \ads \xi \lp  
\mathbb{P}_{\fb}^*\delta\eta\rp \rp.
\end{eqnarray*}
Here we used property~\textbf{(SM3)} of the slice mapping $\Psi$.

Since $\delta \gamma \mapsto i_{\fq}^* (\ads {\delta \gamma} \mu)$  
is an 
isomorphism between $\fq$ and $\fp^*$, we conclude that the partial  
derivative
$D_{\fq}F_1(\zero)$ is an isomorphism.
Thus the Implicit Function Theorem implies that there is a function
$\gamma: \mathcal{U}_1 \subset \mathcal{U} \times \fa \times \fb  
\to \fq$ such
that
\[
F_1(\eta, v, \alpha, \beta, \gamma(\eta, v, \alpha, \beta))
= i_{\fp}^* 
\ads {\xi + \alpha + \beta + \gamma(\eta, v, \alpha, \beta)}  
\cJ(\eta, v) 
= 0
\]
for all $(\eta, v, \alpha, \beta) \in \mathcal{U}_1$. In other  
words, we have 
found a $\fb^*\times V\times\fa\times\fb$--parameter family
of points that satisfy  part $\textbf{(RE1)}$ of the
split critical point equations. Set
\begin{equation}
\label{omega 1}
\omega_1(\eta, v, \alpha, \beta)
:= \xi + \alpha + \beta + \gamma(\eta, v, \alpha, \beta).
\end{equation}

\medskip

\noindent\textbf{Step 2.} In this step we assume that the
subspace $\fb$ is reflexive, that is, $\fb^{**}\simeq\fb$.
We now construct a $\fb^*\times V\times\fa$--parameter family
of points satisfying the relative equilibrium equations
$\textbf{(RE1)}$ and $\textbf{(RE3)}$ by applying the Implicit 
Function Theorem to $\textbf{(RE3)}$, solving for
the $\fb$ component of the family of points constructed in Step 1.

Let $F_2:\mathcal{U}_1\subset\fb^*\times
V\times\fa\times\fb\rightarrow\fb^{**}\simeq\fb$ be
the mapping defined by $F_2(\eta,v,\alpha,\beta):=
D_{\fb^*}\mathcal{F}(\eta,v,\omega_1(\eta,v,\alpha,\beta))$. Since  
we intend 
to solve the equation $F_2 = 0$ for the $\fb$ parameter using the  
Implicit
Function Theorem, we compute $D_{\fb}F_2(0,0,0,0)$.
Given arbitrary $\delta\beta\in\fb$ and $\delta\eta\in\fb^*$,
\begin{eqnarray*}
\langle D_{\fb}F_2(\mathbf{0})\delta\beta,\delta\eta\rangle
&=& \left. \smallfrac {d}{dt} \left. \smallfrac{d}{ds}
\lp \cH - \cJ^{\omega_1(0,0,0,t \, \delta\beta)} \rp(s \, \delta\eta,0)
\right|_{t=0}\right|_{s=0} \\
&=& D \cJ(\zero)(\delta\eta,  
0)(D_{\fb}\omega_1(\mathbf{0})\delta\beta)\\
&=& \langle\mathbb{P}_{\fb}^*\delta\eta,
D_{\fb}\omega_1(\mathbf{0})\delta\beta\rangle \\
&=& \langle \delta \eta, \mathbb{P}_{\fb} (\delta \beta 
+ D_{\fb}\gamma(\mathbf{0})\delta\beta) \rangle \\
&=& \langle \delta \eta, \delta \beta \rangle
\end{eqnarray*}
follows from property \textbf{(SM3)} of the slice map and 
the formula~(\ref{omega 1}) for the generator $\omega_1$.
Hence $D_{\fb}F_2(\mathbf{0})$ is
the identity map. The Implicit Function Theorem thus implies that  
there is
a function $\beta: \mathcal{U}_2 \subset \mathcal{U} \times \fa  
\rightarrow \fb$
satisfying
$F_2(\eta, v, \alpha, \beta(\eta, v, \alpha)) =
D_{\fb^*}\mathcal{F}(\eta,v,\omega_1(\eta,v,
\alpha,\beta(\eta, v, \alpha)))=0$ for all
$(\eta, v, \alpha) \in  \mathcal{U}_2$.
Set
\begin{equation}
\label{omega 2}
\omega_2(\eta, v, \alpha) := \omega_1(\eta, v, \alpha,
\beta(\eta, v, \alpha)).
\end{equation}

\medskip
\noindent\textbf{Step 3.} We now treat the~(\textbf{RE4}) component  
of the 
relative equilibrium equation. We use the standard  
Lyapunov--Schmidt reduction 
procedure of bifurcation theory to partially solve (\textbf{RE4}).

Let $\mathcal{L}: V \to V^*$ denote the linear transformation satisfying
\[
\langle \mathcal{L} \, v, w \rangle := D_{VV} \lp \hxi \rp(\zero)(v, w)
= D^2 \lp \hxi \rp(\zero)((0, v), (0, w))
\]
for all $v$ and $w \in V$. 
Set $V_0 := \ker \, \mathcal{L}$ and choose
closed subspaces $V_1 \subset V$ and $V_2 \subset V^*$ such that
\[
V = V_0 \oplus V_1
\qquad\text{and}\qquad
V^* = {\rm range} \, \mathcal{L} \oplus V_2.
\]

To guarantee the existence of the complements $V_1$ and $V_2$
in infinite dimensions, we assume that $\mathcal{L}$ is a Fredholm  
operator.
Let $\mathbb{P}: V^* \rightarrow V_2$ denote the
 projection determined by the decomposition
of $V^*$. Note that since the operator $\mathcal{L}$ is such that  
for all $v$ and $w \in V$, $\langle \mathcal{L} \, v, w  
\rangle=\langle \mathcal{L} \, w, v \rangle$, the spaces $V_0$ and  
$V_2$ can be naturally identified by choosing an inner product (when  
$V$ is infinite dimensional we still can do it provided that $V$ is  
a Hilbert space).

Define $F_3: \fb^* \times V_0 \times V_1
\times \fa \rightarrow {\rm range} \, \mathcal{L}$
by
\[
F_3(\eta, v_0, v_1, \alpha) := (\mathbb{I} -  \mathbb{P})
D_V \lp \cH - \cJ^{\omega_2(\eta, v_0  +  v_1, \alpha)} \rp 
(\eta, v_0  +  v_1).
\]
Using the Implicit Function Theorem once more, we can solve the  
equation 
$F_3(\eta, v_0, v_1, \alpha) = 0$ for $v_1$. For any $\delta v_1\in  
V_1$,
\[
D_{V_1}F_3(\mathbf{0})\cdot\delta v_1
= (\mathbb{I} -  \mathbb{P}) \lp \mathcal{L} \, \delta v_1
        - D \cJ^{D \omega_2(\zero) (0, \delta v_1, 0)}(\zero) \rp 
= \mathcal{L} \, \delta v_1,
\]
since $(\mathbb{I} -  \mathbb{P}) \mathcal{L} = \mathcal{L}$ and
$V \subset \mbox{ker} \, D \cJ(\zero)$.
Thus $D_{V_1}F_3(\zero)$ is an isomorphism of $V_1$ onto 
${\rm range} \, \mathcal{L}$ and the Implicit Function Theorem  
guarantees
the existence of a neighborhood $\mathcal{U}_3$ of $(0,0,0)
\in \fb^* \times V_0 \times \fa$, and a local function
$v_1:\mathcal{U}_3\rightarrow V_1$
such that
\[
F_3(\eta, v_0, v_1(\eta,v_0,\alpha), \alpha)=0,
\]
for any $(\eta,v_0,\alpha)\in\mathcal{U}_3$.

Define the {\bfi generator map} $\Xi: \mathcal{U}_3 \rightarrow \fg$,
$B: \mathcal{U}_3 \rightarrow V_2$, and
$\rho: \mathcal{U}_3 \rightarrow \fm$ by
\begin{eqnarray*}
\Xi(\eta, v_0, \alpha)
&:=& \omega_2(\eta, v_0  +  v_1(\eta, v_0, \alpha), \alpha) \\
B(\eta, v_0, \alpha)
&:=& \mathbb{P} D_V \lp \cH - \cJ^{\Xi(\eta, v_0, \alpha)} \rp
(\eta, v_0 + v_1 (\eta, v_0, \alpha)) \\
\rho(\eta, v_0, \alpha)
&:=& \iota_{\fm}^* 
\ads {\Xi(\eta, v_0, \alpha))}\cJ ((\eta, v_0  +  v_1(\eta, v_0,  
\alpha)).
\end{eqnarray*}

In a sufficiently small neighborhood $\mathcal{U}_3$ 
of the origin any solution $(\eta, v_0, \alpha)$ of the equations
\begin{equation}
\left\{
\begin{array}{lll}
\text{\textbf{(B1)}}\ \ \ B(\eta, v_0, \alpha) &=& 0,\\
 & & \\
\text{\textbf{(B2)}}\ \ \ \rho(\eta, v_0, \alpha) &=& 0
\end{array}
\right.
\end{equation}
determines a relative equilibrium $\Psi(\eta, v_0  +  v_1(\eta,  
v_0, \alpha))$
with generator $\Xi(\eta, v_0, \alpha)$. On the other hand, any  
relative 
equilibrium $m$ sufficiently near $m_e$ in the slice $\Psi(\fb^*  
\times V)$
satisfies $m = \Psi(\eta, v_0  +  v_1(\eta, v_0, \alpha))$ for some  
solution
$(\eta, v_0, \alpha)$ of (\textbf{B1}) and~(\textbf{B2}); any  
generator $\xi$
of $m$ satisfies $\xi - \Xi(\eta, v_0, \alpha) \in \fg_m$.
Equations~(\textbf{B1}) and~(\textbf{B2}) will be
usually referred to as the {\bfi bifurcation equation} and
the {\bfi rigid residual equation} respectively. Let   $R:\fh\times  
\fm\sus\times V_0\rightarrow \fm\sus\times V_0$ be the mapping that  
groups both equations, that is,
\begin{equation*}
\begin{array}{cccc}
R:&\fh\times \fm\sus\times V_0&\longrightarrow &\fm\sus\times V_0\\
	&(\alpha,\eta,v_0)&\longmapsto&(\rho(\eta, v_0,  
\alpha),B(\eta, v_0, \alpha)).
\end{array}
\end{equation*}
We will refer to the equality 
\begin{equation}
\label{reduced version}
R(\alpha,\eta,v_0)=0
\end{equation}
as the {\bfi reduced critical point equations}.

\begin{remark} \label{rem_gradient}
\normalfont
Notice that even though the critical point  
equations~(\ref{critical}) 
determining the relative equilibria in our situation can be 
naturally understood as a gradient equation, this analytic feature  
is not 
in general available for the reduced version~(\ref{reduced  
version}) of these 
equations.

A particular case where the gradient character of~(\ref{critical})  
is preserved by the reduction procedure is when the relative  
equilibrium $m_e$ that we start with has total isotropy, that is, it  
is actually an equilibrium and $G_{m_e}=G$. Notice that in this  
case  $\fm=\fp=\{0\}$ since $\fg_{m_e} = \fg$. Therefore, the rigid  
residual equation~\textbf{(B2)} is  
trivial and then, as we will show, the bifurcation  
equation~\textbf{(B1)} is a gradient equation. The strategy that we  
will take follows very closely the one introduced  
in~\cite{constrained paper}.

If $\fm=\fp=\{0\}$, then any coordinate chart 
$\psi: {\mathcal U} \subset X \to M$ such that $\psi(0) = m_e$ is a  
slice 
mapping at $m_e$, with $V = X$, and the critical point 
equations~\textbf{(RE1)}--\textbf{(RE4)} collapse to the single equation
\begin{equation}
\label{single equation}
D_V \lp \hxi \rp(v) = 0.
\end{equation}
In this situation only the third step of the general procedure, 
the Lyapunov-Schmidt reduction, is nontrivial.  

We fix an inner product $\langle\cdot,\cdot\rangle$ on $V$ and denote by
$\nabla_V \lp \hxi \rp(v)$ the usual gradient of $\cH - \cJ^\xi$  
with respect to
$\langle\cdot,\cdot\rangle$, i.e. 
\[
\langle\nabla_V \lp \hxi \rp(v),w\rangle = D_V\lp \hxi \rp(v)\cdot w
\]
for any $w \in V$. If $m_e$ is a relative equilibrium with
generator $\xi$, the relative equilibria around $m_e$ are
given by the zeroes of the map $F:V\times\fg\rightarrow V$ defined by
\[F(v,\alpha)=\nabla_V\lp \mathcal{H}-\mathbf{j}^{\xi+\alpha} \rp(v).\]
Let $L:V\rightarrow V$ be the mapping
defined by $L(v)=D_V F(0,0)\cdot v$. It can easily be verified that
\[
\langle L(v),w\rangle=D_{VV}\lp \hxi \rp(0)(v,w)
\]
for any $v$ and $w \in V$.
Note that the mapping $L$ is a self-adjoint operator; hence if we set 
$V_0=\ker L$ and $V_1= \mbox{range} \, L$, then $V$ has the orthogonal 
decomposition $V=V_0\oplus V_1$. Let $\mathbb{P}:V\rightarrow V_0$
denote the canonical projection with respect to the splitting  
$V=V_0\oplus
V_1$. Now, if
we decompose $v\in V$ as $v=v_0+v_1$, with $v_0\in V_0$ and $v_1\in  
V_1$, 
and apply the Implicit Function Theorem to the equation
\[
(\mathbb{I}-\mathbb{P})F(v_0+v_1,\alpha)=0, 
\]
we obtain a function $v_1:V_0\times\fg\rightarrow V_1$ such that
\begin{equation}
\label{definition of v 1}
(\mathbb{I}-\mathbb{P})F(v_0+v_1(v_0,\alpha),\alpha)=0.
\end{equation}
The remaining equation that is, the bifurcation equation, is
\[
B(v_0,\alpha):=\mathbb{P}F(v_0+v_1(v_0,\alpha),\alpha)
=0.
\]
We now show that the map $B$ is the gradient of
$g(v_0,\alpha):= \lp \mathcal{H}-\mathbf{j}^{\xi+\alpha}  
\rp(v_0+v_1(v_0,\alpha))$, that is
\[B(v_0,\alpha)=\nabla_{V_0}g(v_0,\alpha).\]
Indeed, note that for any $w\in V_0$
\begin{eqnarray*}
\langle\nabla_{V_0}g(v_0,\alpha),w\rangle
&=& D_V \lp \cH - \cJ^{\xi + \alpha} \rp (v_0+v_1(v_0,\alpha))\cdot
        (w+D_{V_0}v_1(v_0,\alpha)\cdot w)\\
&=& \langle F(v_0+v_1(v_0,\alpha),\alpha), {\mathbb P} w + 
({\mathbb I} - {\mathbb P}) D_{V_0}v_1(v_0,\alpha)\cdot w) \rangle \\
&=& D_V \lp \cH - \cJ^{\xi + \alpha} \rp (v_0+v_1(v_0,\alpha))\cdot w \\
&=& \langle {\mathbb  P}F(v_0+v_1(v_0,\alpha),\alpha),w\rangle
=\langle B(v_0,\alpha),w\rangle, 
\end{eqnarray*}
since $w \in V_0 = \mbox{range} \, {\mathbb P}$, 
$D_{V_0}v_1(v_0,\alpha)\cdot w \in V_1 = \mbox{range} \, ({\mathbb  
I} - {\mathbb P})$,
$\mathbb{P}$ is self-adjoint, and (\ref{definition of v 1}) is  
satisfied. \ \ \ $\blacklozenge$

\end{remark}

\subsection{The equivariance properties of the reduced critical  
point equations}
\label{equivariance bifurcation}

The symmetries of the relevant
equations play an important role the solution of a bifurcation problem
(see for instance~\cite{g2}). We will see that
if the $G$--action on $M$ is {\bfi proper}, then
the relative equilibrium equations~\textbf{(B1)} and~\textbf{(B2)}  
can be 
constructed so as to be equivariant with respect to the induced  
action of 
$\sym$ on $\fb^* \times V_0$. 
Here $G_\xi$ denotes the isotropy subgroup of the generator $\xi\in\fg$
of the relative equilibrium $m_e\in M$ with respect to the adjoint  
action 
of $G$ on $\fg$.

An {\bfi equivariant slice mapping} is a mapping 
$\Psi:\mathcal{U}\subset\fm^*\times V\rightarrow M$ satisfying 
Definition~\ref{slice mapping} and the condition
\begin{description}
\item[(ESM)] 
The subspace $\fm^*$ of $\fg^*$ is ${\rm Ad}^*_{(\sym)}$--invariant and
the slice mapping $\Psi:\mathcal{U}\subset\fm^*\times V\rightarrow M$ is
$\sym$--equivariant with respect to the coadjoint action of 
$\sym$ on $\fm^*$ and some action of $\sym$ on $V$.
\end{description}
Note that since the group $\sym$ is compact
and fixes $(0,0)\in\fm^*\times V$, the open neighborhood
$\mathcal{U}$ of $(0,0)\in\fm^*\times V$ in~\textbf{(ESM)}
can always be chosen to be $\sym$--invariant.

\begin{proposition}
If the group $G$ acts properly on $M$ and the coordinate chart  
$\psi: {\cal U} \subset X \to M$ with $\psi(0) = m_e$ is  
equivariant with 
respect to some action of $\sym$ on $X$, then the subspaces 
$\fb$, $\fq$, $V$, and $W$ can be taken to be $\sym$ invariant. For
these choices, the slice mapping constructed in  
Proposition \ref{build_slice} is $\sym$--equivariant. 
\end{proposition}
\proof
First, we show that $\sym$--invariant decompositions
$\fg=\fg_{m_e} \oplus\fm\oplus\fp$ and $X = V \oplus W$ exist.
Note that the isotropy subgroup $G_{m_e}$ is compact,
since the action of $G$ on $M$ is assumed to be proper; 
consequently the subgroup $\sym$ is also compact. This 
guarantees the existence of a ${\rm Ad}_{(\sym)}$--invariant inner
product on $\fg$, which we can use to determine a 
${\rm Ad}_{(\sym)}$--invariant decomposition
$\fg=\fg_{m_e} \oplus\fm\oplus\fp$ of the Lie algebra. 

The orthogonal complement to $\fg_\mu\cdot m_e$ in
$\ker T_{m_e}\mathbf{J}$ with respect to a $\sym$--invariant
inner product is an invariant subspace. Hence the preimage with  
respect to 
the equivariant map $T_0 \psi$ of this orthogonal complement is a 
$\sym$--invariant subspace of $X$; we choose this subspace as
the vector space $V$ in Definition~\ref{slice mapping}. Analogously, the 
space $W$ can
be chosen to be invariant under the action on $X$.

Given these choices of subspaces, the action of $\sym$ on 
$M$ induces a well-defined action on $\fm^* \times V$ via the slice map.
Equivariance of the momentum map, the coordinate chart, and the  
projection 
${\mathbb P}_\fm$ imply that the slice map $\Psi$ is equivariant.
\qed

Recall that the relative equilibrium equations were obtained using two 
consecutive applications of the Implicit Function Theorem (Steps 1  
and 2) 
followed by the Lyapunov-Schmidt reduction procedure (Step 3). 
It is well known that if the Implicit Function Theorem is applied to an
equation $F = c$ determined by an equivariant map $F$ and a fixed  
point $c$ of
the group action, then the resulting implicitly defined function is  
also 
equivariant. In addition, if the Lyapunov-Schmidt reduction 
procedure is applied to such an equation using invariant subspaces,
then the resulting functions and equations will be equivariant. 
(See, e.g., \cite{g1, g2} for precise statements and proofs of  
these results).  Using these fundamental results, we now show that,
given appropriate choices of slice maps and subspaces, the  
generator map $\Xi$
and the functions $B$ and $\rho$ determining the reduced relative  
equilibrium 
equations are equivariant with respect to the induced $\sym$ action on 
$\fb^* \times V \times \fg_{m_e}$.

\begin{proposition}
\label{equivariance features}
If the spaces $\fb$, $\fq$, $V$, and $W$ are $\sym$ invariant and
the slice mapping is $\sym$--equivariant, then the 
maps $\Xi$, $v_1$, $B$,  $\rho$, and $F$ are all $\sym$--equivariant.
\end{proposition}
\proof
It suffices to show that the functions $F_1$, $F_2$,
and $F_3$ given in steps 1, 2, and 3 are $G_{m_e}\cap  
G_\xi$--equivariant.
We first consider the mapping 
$F_1: \mathcal{U} \times \fg_{m_e} \times \fm \times \fp \to \fp^*$
introduced in Step 1. For arbitrary $g\in G_{m_e}\cap G_\xi$:
\begin{eqnarray*}
F_1(g\cdot (\eta, v, \alpha, \beta, \gamma))
&=& i_\fp^* \, \ads {\xi + g\cdot \alpha + g\cdot\beta + g \cdot  
\gamma)} 
\cJ(g\cdot\eta,g\cdot v) \\
&=& i_\fp^* \, \ads {\Ad g (\xi + \alpha + \beta + \gamma)} 
\Ads {g^{-1}} \cJ(\eta, v) \\
&=& i_\fp^* \, \Ads {g^{-1}} \lp
\ads {\xi + \alpha + \beta + \gamma} \cJ(\eta, v) \rp \\
&=& \Ads {g^{-1}} \lp i_\fp^* \, 
\ads {\xi + \alpha + \beta + \gamma} \cJ(\eta, v) \rp \\
&=& g \cdot F_1(\eta, v, \alpha, \beta, \gamma).
\end{eqnarray*}
Thus $F_1$ is $G_{m_e}\cap G_\xi$--equivariant and, hence, the
implicitly defined functions $\gamma$ and $\omega_1$ are also
$G_{m_e}\cap G_\xi$--equivariant. An analogous verification
can be carried out for the mapping $F_2$ in Step 2, allowing
us to conclude that the function $\omega_2$ is also $G_{m_e}\cap
G_\xi$--equivariant.

To establish the invariance (respectively equivariance) of the  
spaces and 
maps constructed in Step 3, we first note that $\hxi$ is 
$\sym$--invariant, since the augmented Hamiltonian $h - J^\xi$ 
is $G_\xi$--invariant and the slice map $\Psi$ is 
$\sym$--equivariant. Equivariance of the map $F$, 
and hence invariance of the subspaces ${\rm ker} \, F$ and 
${\rm range} \, F$, follows immediately from the invariance of 
$\hxi$. The compactness of the group $\sym$ allows us to choose
$\sym$--invariant complements $V_1$ and $V_2$ to
${\rm ker} \, F$ and ${\rm range} \, F$.
(See for instance~\cite[Proposition 2.1]{g2}.) 
With these choices, the canonical projection $\mathbb{P}$
and the function $F_3$ are equivariant. Consequently the function
$v_1$, as well as the generator map $\Xi$ and the reduced relative  
equilibrium
equations are equivariant, as required.
\qed

\subsection{Treatment of the rigid residual equation}
\label{Slice map refinements and the rigid residual equation}

In this section we consider some situations in which the rigid
residual map is either trivial or at least fairly simple.
For example, if $G$ is Abelian, then the full rigid equation
$\ads \xi \J(m_e) = 0$ is trivial. Hence, the rigid residual
equation is obviously satisfied as well.
If $G$ is not Abelian, an appropriate choice
of a slice map $\Psi: \fb^* \times V\rightarrow M$ can be useful in
identifying the solutions of the residual rigid equation.
We will present  a few cases in which these helpful
choices are possible.

Given a relative equilibrium $m_{e}$ of $h: M \times
\fg \rightarrow \mathbb{R}$, let $\mu := \J(m_{e}) \in \fg^*$.
Let $\mathcal{O}_{\mu}\subset\fg^{*}$ be the
coadjoint orbit through $\mu\in\fg^*$, with
tangent space $T_{\mu}\mathcal{O}_{\mu}$ at $\mu$ given by
\[T_{\mu}\mathcal{O}_{\mu}=\{{\rm ad}^*_{\zeta}
\mu\mid \zeta\in\fg\}.\]
We shall say that a subspace $\fq\subset
\fg$ is $\fg_{\mu}$--invariant if
$[\fg_{\mu},\fq]\subset\fq$.

We now prove that, generically, the rigid equation $\rho$ can be
reduced by an appropriate choice of slice map to an equation on  
$\fg_{\mu}$.

\begin{proposition}
\label{reduction to smaller}
If the complement $\fq$ to $\fg_{\mu}$ in $\fg$ is  
$\fg_{\mu}$--invariant, 
then given any slice map $\Psi:\mathcal{U} \rightarrow M$ at $m_{e}$, 
there exists a map $\phi:\tilde\mathcal{U} \rightarrow\fq$ such that
\begin{enumerate}
\item
the map $\tilde\Psi:\tilde\mathcal{U}\subset\mathcal{U}\rightarrow M$ 
given by
\begin{equation}
\label{modified slice}
\tilde\Psi(\eta,v)=\exp(\phi(\eta,v))\cdot\Psi(\eta,v)
\end{equation}
is also a slice map
\item
the associated generator map $\tilde\Xi$ takes values in
$\fg_{\mu}$ 
\item
the pullback $\tilde \cJ := \J \circ \tilde \Psi$ of the
momentum map takes values in $\mu + \fq^\circ$
\item
$\phi(0,0)=0$ and $D\phi(0,0)=0$. 
\end{enumerate}
If the original slice mapping is $\sym$--equivariant, then $\tilde\Psi$
is equivariant.
\end{proposition}

\noindent\textbf{Proof} 
We obtain the map $\phi$ through yet another application of the Implicit
Function Theorem. Define $R: \fb^* \times V \times \fq \rightarrow  
\fq^*$ by
\begin{equation}
\label{phi defined}
R(\eta, v, \phi)
= i_\fq^* \left(\bJ(\exp(\phi) \cdot \Psi(\eta, v)) - \mu \right),
\end{equation}
with differential
\begin{eqnarray*}
DR(\zero)(\delta \eta, \delta v, \delta\phi)
&=& i_{\fq}^*D\bJ(m_{e})\lp T_{(0,0)}\Psi(\delta\eta,\delta v)
+ \delta \phi)_M(m_{e}) \rp \\
&=& i_{\fq}^* \lp \mathbb{P}_{\fb}^* \delta\eta - \ads {\delta  
\phi} \mu \rp \\
&=& - i_{\fq}^* \ads {\delta \phi} \mu 
\end{eqnarray*}
for arbitrary $\delta\eta\in\fb^*$, $\delta v\in V$, and $\delta  
\phi \in \fq$.
Here~\textbf{(SM3)}, equivariance of the momentum map, and the identity 
$i_{\fq}^* \mathbb{P}_{\fb}^* = \lp \mathbb{P}_{\fb} \circ i_{\fq}  
\rp^* = 0$
have been used to simplify the expressions.
Since $\eta \mapsto i_{\fq}^* \ads \eta \mu$ is an isomorphism from  
$\fp$ to
$\fq^*$, the Implicit Function Theorem implies
that there is a neighborhood $\tilde \mathcal{U}$ of $(0, 0)$
in $\fb^* \times V$ and a function $\phi: \tilde \mathcal{U} \to  
\fq$ such 
that $\phi(0,0)=0$, $D \phi(0, 0) = 0$, and $R(\eta, v, \phi(\eta,  
v)) = 0$.

Using $\phi:\tilde\mathcal{U}\subset\fb^*\times V \to \fq$ and 
(\ref{modified slice}), we see that the pullback $\tilde \cJ$ 
of the momentum map satisfies
\[
i_{\fp}^* \lp \ads {\xi + \alpha + \beta}\tilde \cJ(\eta,v) \rp
= i_{\fp}^* \lp \ads {\xi + \alpha + \beta}(\tilde \cJ(\eta,v) -  
\mu) \rp
= 0
\]
for all $(\eta,v,\alpha,\beta)\in\tilde \mathcal{U}_1$.
Thus executing Step 1 of Section~\ref{the bifurcation equations}  
using the 
modified slice mapping $\tilde\Psi$ yields a mapping 
$\tilde\gamma:\tilde\mathcal{U}_1\subset
\fb^*\times V\times\fa\times\fb\rightarrow\fq$ satisfying 
\begin{eqnarray*}
0 &=& F_1(\eta,v,\alpha,\beta,\tilde\gamma(\eta, v, \alpha,\beta))\\
&=& i_{\fp}^* \lp \ads {\xi + \alpha + \beta + \tilde\gamma(\eta,v,
\alpha,\beta)}\tilde \cJ(\eta,v) \rp \\
&=& i_{\fp}^* \lp 
        \ads {\tilde \gamma(\eta,v, \alpha,\beta)}\tilde \cJ(\eta,v) \rp
\end{eqnarray*}
for any $(\eta,v,\alpha,\beta)\in\tilde \mathcal{U}_1$. 
$\tilde\gamma \equiv 0$ clearly satisfies this equation; hence it  
is the 
unique solution of the equation $F_1 \equiv 0$ given by the  
Implicit Function 
Theorem. Thus steps 2 and 3 yield the generator map
\[
\tilde\Xi(\eta,v_0,\alpha) 
= \xi + \alpha + \beta(\eta,v_0 + v_1(\eta,v_0,\alpha),\alpha) \in  
\fg_\mu.
\]

Suppose now that the slice map $\Psi$ has the  
property~\textbf{(ESM)}.  
Note that for any $(\eta,v,\phi)\in \fm^*\times V\times\fq$ and any 
$h\in \sym\subset G_\mu$
\begin{eqnarray*}
R(h\cdot\eta,h\cdot v,h\cdot \phi)
&=& i_\fq^* \left(\bJ(\exp(h\cdot\phi) \cdot \Psi(h\cdot\eta,  
h\cdot v)) 
	- \mu \right)\\
&=& h\cdot i_\fq^* \left(\bJ(\exp(\phi) \cdot \Psi(\eta, v)) 
	- \mu \right) \\
&=& h\cdot R(\eta,v,\phi).
\end{eqnarray*}
Equivariance of $R$ implies that $\phi$, and hence $\tilde \Psi$, are
equivariant.
$\blacksquare$

\medskip

If the hypotheses of Proposition~\ref{reduction to smaller} are  
satisfied, 
the rigid residual equation involves only elements of $\fg_\mu$ and  
$\fg_\mu^*$.
Specifically, if we let $[\ ,\ ]_\mu$ denote the Lie bracket on  
$\fg_\mu$
and $\bJ_\mu: M\to \fg_\mu^*$ denote the momentum map associated to  
the action 
of $G_\mu$ on $M$, namely $\bJ_\mu = i_{\fg_\mu}^* \bJ$, then  
$\rho$ satisfies
\begin{equation}
\label{all goes to smaller}
\rho(\eta, v_0, \alpha) \cdot \beta
= \bJ_\mu(\tilde \Psi(\eta, v_0 + v_1(\eta,v_0,\alpha))) \cdot
[ \tilde \Xi(\eta, v_0, \alpha), \beta ]_\mu,
\end{equation}
for all $\beta \in \fb$. In particular, if $\fg_\mu$ is Abelian, then
$\rho$ is identically zero. Thus we have established the following  
corollary.

\begin{corollary}
\label{abelian matters}
Let $m_e$ be a relative equilibrium with momentum $\mu = \bJ(m_e)$. If 
$\fg_\mu$ is Abelian and there exists a $\fg_\mu$--invariant  
complement to 
$\fg_\mu$ in $\fg$, then there is a slice map with respect to 
which the rigid residual map $\rho$ is identically zero.
\end{corollary}

\medskip

Another approach to the search for solutions of the rigid residual  
equation
is to restrict this search to fixed point subspaces corresponding  
to subgroups 
of the symmetry group of $\rho$. More explicitly, suppose that the  
hypotheses of
Proposition~\ref{reduction to smaller} are satisfied  and that we  
start with an 
equivariant slice map $\Psi$. In that case, 
Proposition~\ref{equivariance features} guarantees that $\rho$ is 
$\sym$--equivariant and satisfies~(\ref{all goes to smaller}).  
Equivariance 
implies that for any Lie subgroup $K\subset\sym$, the map $\rho$  
maps the
set of fixed points of $K$ into the set of fixed points of $K$ in  
$\fm^*$.
Hence all zeroes of the restriction 
\[
\rho^K:(\fm\sus)^K\times V_0^K\times\fg_{m_e}^K\longrightarrow  
(\fm\sus)^K,
\]
of $\rho$ to $(\fm\sus)^K\times V_0^K\times\fg_{m_e}^K$ are also  
zeroes of $\rho$
(here the superscript $K$ denotes the subspace of $K$--fixed points  
with 
respect to the relevant action).  In other words, we can look for   
the solutions of the rigid residual equation by searching the zeroes  
of  its restrictions to different sets of $K$--fixed points, with  
$K$ and arbitrary subgroup of $\sym$ which, in principle, should be  
easier, since the dimension of the system has been lowered without 
introducing additional complexity into the equations.

If the restriction of the Lie bracket of the 
Lie algebra $\g_\mu$ to $\g_\mu^K$ is trivial, then the entire subspace 
$(\fm\sus)^K\times V_0^K\times\fg_{m_e}^K$ consists of solutions of  
the rigid 
residual equation. Indeed, for any $(\eta, v_0, \alpha)\in 
(\fm\sus)^K\times V_0^K\times\fg_{m_e}^K$, if we let 
\[
\nu = \bJ_\mu(\tilde \Psi(\eta, v_0 + v_1(\eta,v_0,\alpha))) 
\qquad \mbox{and} \qquad
\zeta = \tilde \Xi(\eta, v_0, \alpha),
\]
then 
\[
\rho(\eta, v_0, \alpha) 
= \langle\nu,  [\zeta, \cdot ]_{\fg_\mu}\rangle.
\]
The equivariance of $\tilde \Xi$ and $\bJ_\mu$ implies that
$\zeta \in \fg_\mu^K$ and $\nu \in (\fg_\mu^*)^K$. 
Also, since $\fm\subset\fg_\mu$, we have 
$(\fm)^K\subset(\fg_\mu)^K$. Therefore, since $(\fm^*)^K \simeq  
(\fm^K)^*$, we have for any $\xi\in\fm^K$
\[
\langle\rho(\eta, v_0, \alpha),\,\xi\rangle=\langle\nu, [\zeta, \xi  
]_{\fg_\mu^K}\rangle=0,
\]
due to the hypothesis on the Lie bracket on $\fg_\mu^K$.
The arbitrary character of $\xi\in\fm^K$ implies that $\rho(\eta,  
v_0, \alpha)=0$.

Thus we have then proved the following

\begin{proposition}
Let $m_e$ be a relative equilibrium with momentum $\mu = \bJ(m_e)$  
and generator $\xi\in\fg$. If 
there exists a $\fg_\mu$--invariant complement to 
$\fg_\mu$ in $\fg$, then for any subgroup $K\subset\sym$ for which  
the restriction of the Lie bracket of the Lie algebra $\g_\mu$ to  
the set of fixed points $\g_\mu^K$ is trivial, there is a slice map  
$\tilde \Psi$ with respect to 
which the entire subspace $(\fm\sus)^K\times  
V_0^K\times\fg_{m_e}^K$ consists of zeroes of the rigid residual  
equation $\rho$.
\end{proposition}

(See~\cite{rosou} for persistence results on nondegenerate 
Hamiltonian relative equilibria valid under conditions of this sort.)

\section{Persistence in Hamiltonian systems with Abelian symmetries}
\label{Persistence in Hamiltonian systems with Abelian symmetries}

In this section we will focus on the relative equilibria of Hamiltonian
systems with Abelian symmetry groups. Specifically, we consider a  
Hamiltonian
system $(M,\omega,h,G,\mathbf{J}:M\rightarrow\fg^{*})$, where  
$\omega$ is
the symplectic structure of $M$, the Abelian group $G$ acts  
properly on $M$
with associated momentum map ${\mathbf J}$, and the Hamiltonian $h$  
is $G$--invariant. Let $m_e \in M$ be a relative equilibrium with  
generator $\xi$ and
momentum $\mu = \bJ(m_e)$. Since the adjoint and coadjoint actions  
of an 
Abelian group are trivial, $G_\mu=G$ and the rigid residual 
equation~\textbf{(B2)} is trivially satisfied. We also assume that the 
bifurcation equation~\textbf{(B1)} is trivial, i.e.  
that $m_e$ is a {\bfi nondegenerate relative equilibrium}, with
\[
\ker \, D^2 (h - \bJ^\xi)(m_e) = \fg_{\mu} \cdot m_e = \fg \cdot m_e.
\]
In this situation Steps 1 through 3 in Section~\ref{the bifurcation
equations} guarantee the
existence of a $\fm^*\times\fg_{m_e}$--parameter family of relative
equilibria
{\bfi persisting} from $m_e$, whose dimension and structure we will  
study in
what follows.
We use the word {\bfi persistence} as opposed to the word {\bfi  
bifurcation},
given that the
latter is customarily used to indicate a qualitative change in the  
family of 
relative equilibria as a given parameter is varied.  This is  
analytically reflected in the need for a nontrivial Lyapunov-Schmidt  
reduction procedure in order to write the bifurcation equations. We  
shall see that in the 
case at hand no such tool will be necessary. 

In this section we will use a very special slice mapping based on the
Marle--Guillemin--Sternberg normal form~\cite{nfm, nfgs, gs} (we will
refer to it as the {\bfi MGS--normal form}), that we briefly describe.
The following exposition includes without proof the details of
the MGS--normal form that will be needed in our discussion. For
additional information the reader should consult the above  
mentioned original
papers or~\cite{rosou, thesis, redham}.

We start by introducing the main ingredients of the MGS
construction. Even though we are in the Abelian case we will  
present, for future reference,  the general case. First, the  
properness of the $G$--action implies that the isotropy
subgroup $G_{m_e}$ is compact.
Second, the vector space $V_{m_e}:=T_{m_e}(G\cdot  
m_e)^\omega/(T_{m_e}(G\cdot
m_e)^\omega\cap T_{m_e}(G\cdot m_e))=\ker  
T_{m_e}\mathbf{J}/T_{m_e}(G_\mu\cdot
m_e)$ is called the {\bfi symplectic
normal space\/}, which is a
symplectic vector space with the symplectic normal form
$\omega_{V_{m_e}}$ defined by
\[\omega_{V_{m_e}}([v],[w]):=\omega(m_e)(v,w),\]
for any $[v]=\pi(v)$ and $[w]=\pi(w)\in V_{m_e}$, and where
$\pi: \ker T_{m_e}\mathbf{J}\rightarrow \ker  
T_{m_e}\mathbf{J}/T_{m_e}(G_\mu\cdot
m_e)$ is the
canonical projection. Let $H:=G_{m_e}$ be the isotropy subgroup of
$m_e$. The mapping $(h,[v])\longmapsto [h\cdot v]$, with
$h\in H$ and $[v]\in V_{m_e}$,  defines a canonical action of the Lie
group $H$ on $(V_{m_e},\omega_{V_{m_e}})$, where $g\cdot u$ denotes the
tangent lift of the $G$--action on $TM$, for $g \in G$ and $u \in
TM$. The canonical $H$--action on $V_{m_e}$ is linear by  
construction and
 globally Hamiltonian with
momentum map $\mathbf{J}_{V_{m_e}}:V_{m_e}\rightarrow \fh^*$ given by
\[
\langle\mathbf{J}_{V_{m_e}}(v),\eta\rangle
=\frac{1}{2}\omega_{V_{m_e}}(\eta_{V_{m_e}}(v),v),
\]
for arbitrary $\eta\in \fh$ and $v\in V_{m_e}$.Here,  
$\eta_{V_{m_e}}(v) =
\eta \cdot v$ is the induced Lie algebra representation of
$\mathfrak{h}$ on $V_{m_e}$.

The MGS--normal form is based on the construction of a model  
$(Y,\omega_{Y})$
for the symplectic
$G$--manifold $(M,\omega)$ that we introduce in the following  
proposition.

\begin{proposition}
\label{marle tube}
Let $(M,\omega)$ be a symplectic manifold and let $G$ be a
Lie group acting properly on $M$ in a globally Hamiltonian fashion,
with invariant momentum map $\mathbf{J}:M\rightarrow
\fg^*$. Let $m_e\in M$ and denote
$\mathbf{J}(m_e)=\mu\in\fg^*$. Let
$(V_{m_e},\omega_{V_{m_e}})$ be the symplectic normal space at  
$m_e\in M$.
Relative to an  ${\rm Ad}_{G_{m_e}}$--invariant inner product on $\fg$
consider the inclusions $\fm^*\subset \fg_{\mu}^{*}\subset
\fg^{*}$. Then, the manifold
\[Y:=G\times_H(\fm^*\times V_{m_e})\]
can be endowed with a symplectic structure $\omega_{Y}$ with
respect to which the left $G$--action
$g\cdot[h,\eta,v]=[gh,\eta,v]$ on $Y$ is globally
Hamiltonian with  momentum map
$\mathbf{J}_{Y}:Y\rightarrow \fg^{*}$  given by
\begin{equation}
\label{marle momentum}
\mathbf{J}_{Y}([g,\rho,v])= {\rm Ad}_{g^{-1}}\sus(\mu + \rho +  
\mathbf{J}_{V_{m_e}}(v)).
\end{equation}
\end{proposition}

\begin{theorem}[Marle-Guillemin-Sternberg Normal Form]
\label{mnf}
Let $(M,\omega)$ be a symplectic manifold and let $G$ be a Lie
group acting properly on $M$ in a globally Hamiltonian fashion,
with associated invariant momentum map
$\mathbf{J}:M\rightarrow \fg^{*}$. Let $m_e\in M$ and denote
$\mathbf{J}(m_e)=\mu\in\fg^{*}$,
$H:=G_{m_e}$. Then the manifold
\[Y:=G\times_H(\fm^*\times V_{m_e})\]
introduced in Proposition~\ref{marle tube}
is a Hamiltonian $G$--space and there are $G$--invariant
neighborhoods $U$ of $m_e$ in $M$, $U'$ of $[e,0,0]$ in $Y$, and
an equivariant symplectomorphism $\phi:U\rightarrow U'$ satisfying
$\phi(m_e)=[e,0,0]$ and $\mathbf{J}_{Y}\circ\phi=\mathbf{J}$.
\end{theorem}

Since we intend to prove general statements about relative equilibria of
Hamiltonian systems with Abelian symmetries, the previous theorem allows
us to reduce the problem to the study of systems of the form  
$(Y,\omega_Y)$.
Indeed, we will assume that the MGS--normal form is constructed  
around the
relative equilibrium $m_e$ represented by $[e,0,0]$ in ``MGS  
coordinates''. 
It can be easily shown that the map given by
\begin{equation}
\label{marle slice mapping}
\begin{array}{cccc}
\Psi:&\fm^*\times V_{m_e}&\longrightarrow&Y\\
 &(\eta,v)&\longmapsto&[e,\eta,v]
\end{array}
\end{equation}
is a slice mapping at the point $[e,0,0]$ for
$D(h_Y-\mathbf{J}_Y^\xi)$,
which is the representation given by the MGS--normal form of
$D(h-\mathbf{J}^\xi)$.

Before stating the following theorem, we recall 
from elementary differential geometry the basic notion 
of the {\bfi rank of a surface} given in parametric form. Let
$g:U\subset\mathbb{R}^n
\rightarrow\mathbb{R}^m$ be a parameterization of a surface $S$
in $\mathbb{R}^m$.
Given a value $u\in\mathbb{R}^n$ of the parameter, the rank of the  
surface
$S_{g(u)}$ at the
point $g(u)\in\mathbb{R}^m$ is the rank of the Jacobian of the  
function $g$ at
$u$. If this
rank is constant, the Fibration Theorem~\cite[Theorem 3.5.18]{mta}  
guarantees
that $S$ is a
submanifold of $\mathbb{R}^m$ and its rank coincides with the  
dimension of $S$
as a manifold
on its own.

\begin{theorem}
\label{abelian bifurcation theorem}
Let $(M,\omega,h)$ be a Hamiltonian system, and $G$ be an Abelian Lie
group acting properly on
$M$ in a globally Hamiltonian fashion, with invariant momentum map 
$\mathbf{J}:M\rightarrow\fg^{*}$.
Suppose that $h$ is $G$--invariant and that the Hamiltonian system  
determined
by $h$ has a nondegenerate relative equilibrium at the point
$m_e\in M$, with generator $\xi\in\fg$. Set $H:=G_{m_e}$ and
$\mu = \mathbf{J}(m_e)\in\fg^{*}$. Then there is a surface $S$ of
relative equilibria through $m_e$ that can be locally expressed as
\begin{equation*}
S=\{[g,\eta,v(\eta,\alpha)]\in Y\mid g\in G,\eta\in\fm^*,\alpha\in\fh\},
\end{equation*}
using the MGS normal form $Y$ constructed around the orbit $G\cdot  
m_e$. Here
$v:\fm^*\times\fh\rightarrow V_{m_e}$ is a smooth function such  
that $v(0,0)=0$ 
and ${\rm rank}(Dv(\eta,\alpha))= \dim H- \dim H_{v(\eta,\alpha)}$.  
The rank,
${\rm rank} \, S_{[g,\eta,v(\eta,\alpha)]}$, of the surface $S$ at the
relative equilibrium $[g,\eta,v(\eta,\alpha)]$ equals
\begin{equation}
\label{dimension formula for relative}
{\rm rank} \, S_{[g,\eta,v(\eta,\alpha)]}=2(\dim G- \dim H) + (\dim  
H- \dim
H_{v(\eta,\alpha)}).
\end{equation}
\end{theorem}

\noindent\textbf{Proof}  The surface $S$ of relative
equilibria is constructed in Steps 1 through 3 of
Section~\ref{the bifurcation equations}, taking as slice mapping  
the map 
$\Psi(\eta,v)=[e,\eta,v]$ constructed with the help of the  
MGS--normal form.
Indeed, since the nondegeneracy of $m_e$ and the
Abelian character of $G$ imply that~\textbf{(B1)} and~\textbf{(B2)} are
trivially satisfied, there is a neighborhood $\mathcal{U}\subset
\fm^*\times \fh$ of the point $(0,0)$ and functions
$v:\mathcal{U}\rightarrow V_{m_e}$ and $\Xi:\mathcal{U}\rightarrow
\fg$ such that for any $(\eta,\alpha)\in\mathcal{U}$, the point
$[e,\eta,v(\eta,\alpha)]\in Y\simeq M$ is a relative equilibrium of  
the system
$(M,\omega,h)$ with generator $\Xi(\eta,\alpha)\in\fg$.
At the same time, since the Lie group $G$ is Abelian and the  
Hamiltonian flow
$F_t$ associated to $h$ is $G$--equivariant, it is easy to verify that
if the point $[e,\eta,v(\eta,\alpha)]$ is a relative equilibrium
with generator $\Xi(\eta,\alpha)\in\fg$ then,
for any $g\in G$, the point $[g,\eta,v(\eta,\alpha)]$ is
also a relative equilibrium with the same generator.
In order to prove~(\ref{dimension formula for relative}),
we compute $Dv(\eta,\alpha)$ by implicit differentiation
of the equation $F_3(\eta,v(\eta,\alpha),\alpha)=0$ defining the  
function $v$
in Step 3. Note that in this case the space $V_0$ is trivial and we have
dropped the subscript from $v_1$.
Note that $\fq$ is trivial in the Abelian case and hence 
\[
\omega_2(\eta, v, \alpha) = 
\omega_1(\eta, v, \alpha, \beta(\eta, v, \alpha)) 
= \xi + \alpha + \beta(\eta, v, \alpha).
\]
For $u\in V_{m_e}$, for arbitrary $\delta\alpha\in\fh$, if we set 
$\alpha_t = \alpha + t \, \delta \alpha$, we have
\begin{align}
0&= \lp D_{V_{m_e}}F_3(\eta,v(\eta,\alpha),\alpha)\cdot
	(D v(\eta, \alpha)(0 ,\delta\alpha)) \rp \cdot u\notag\\
 &=\ddto D\left(\cH-\cJ^{\xi + \alpha_t 
\beta(\eta, v(\eta, \alpha_t), \alpha_t)}\right)(\eta, 
v(\eta, \alpha_t))\cdot(0,u)\notag\\
 &=D^2\left(\cH-\cJ^{\xi + \alpha +
\beta(\eta, v(\eta, \alpha), \alpha)}\right)
(\eta,v(\eta,\alpha))\cdot((0, Dv(\eta,\alpha)
	\cdot(0, \delta\alpha)), (0, u)) \notag \\
 &\ \ \ - D\mathbf{J}_{V_{m_e}}^{\delta \alpha}(\eta,\alpha))\cdot u.
\label{master equation}
\end{align}
The last equality follows from the identity
\[
\mathbf{j}(\eta, v) = \mu +  \eta +  \cJ_{V_{m_e}}(v),
\]
which implies that
\[
D \mathbf{j}^{\delta \alpha + \delta \beta}(\eta, v) \cdot (0, u)
= \langle D\mathbf{J}_{V_{m_e}}(v))\cdot u, \delta \alpha + \delta  
\beta \rangle
= \langle D\mathbf{J}_{V_{m_e}}(v))\cdot u, \delta \alpha \rangle.
\]
By hypothesis, the quadratic form
\begin{equation}
\label{non degenerate thing we have}
D^2\left(\cH - \cJ^{\xi}\right)(0,0)|_{(\{0\}\times
V_{m_e})\times(\{0\}\times V_{m_e})}
\end{equation}
is nondegenerate; therefore, since nondegeneracy is an open condition,
\begin{equation}
\label{non degenerate thing we want}
D^2\left(\cH - \cJ^{\xi + \alpha + \beta(\eta,\alpha)}
\right)(\eta,v(\eta,\alpha))|_{(\{0\}\times
V_{m_e})\times(\{0\}\times V_{m_e})}
\end{equation}
is nondegenerate
for any $(\eta,\alpha)\in\fm^*\times\fh$ sufficiently close to $(0,0)$.
Hence the rank of $D_{\fh}v(\eta,\alpha)$ equals the 
rank of $D\mathbf{J}_{V_{m_e}}(v(\eta,\alpha))$ at a point 
$(\eta,\alpha)\in\fm^* \times\fh$ sufficiently close to $(0,0)$. Thus 
\begin{equation}
{\rm rank}(D_{\fh}v(\eta,\alpha))={\rm rank}
(D\mathbf{J}_{V_{m_e}}(v(\eta,\alpha)))=\dim
(\fh_{v(\eta,\alpha)})^{{\rm ann}(\fh^*)}=\dim H- \dim  
H_{v(\eta,\alpha)},
\label{rank_equation}
\end{equation}
as required.

The expression (\ref{dimension formula for relative}) for the rank  
of the 
surface $S$ at a relative equilibrium $[g,\eta,v(\eta,\alpha)]$ is a 
straightforward consequence of the formula \ref{rank_equation}
for the rank of $D_{\fh}v(\eta,\alpha)$. The rank of $S$ 
at $[g,\eta,v(\eta,\alpha)]$ is the rank of the parameterization
\[\begin{array}{cccccc}
\mathcal{S}:&G\times\fm^*\times\fh&\longrightarrow&G\times\fm^*
\times V_{m_e}&\longrightarrow&G\times_H(\fm^*\times V_{m_e})\\
  &(g,\eta,\alpha)&\longmapsto&(g,\eta,v(\eta,\alpha))&\longmapsto&[g,\eta,v(\eta,\alpha)]
\end{array}\]
of the surface $S$. The map $\mathcal{S}$ has rank 
\begin{align*}
{\rm rank}(T_{(g,\alpha,\eta)}\mathcal{S})&={\rm
rank}(S_{[g,\eta,v(\eta,\alpha)]})\\
&=\dim G + \dim \fm^* + {\rm rank}(D v(\alpha))- \dim H\\
&=2(\dim G- \dim H) + \dim H- \dim H_{v(\eta,\alpha)},
\end{align*}
at $[g,\eta,v(\eta,\alpha)]$, as required.\ \ \ $\blacksquare$
\medskip

As a corollary to the previous theorem, we can
formulate a generalization of a result due to E. Lerman and S.
Singer~\cite{singreleq},
originally stated for toral actions, to proper actions of Abelian  
Lie groups.
This result was already presented in~\cite{thesis}.

\begin{corollary}
\label{symplectic manifold of equilibria}
Under the hypotheses of Theorem~\ref{abelian bifurcation theorem},  
there is a 
symplectic manifold $\Sigma$ of relative equilibria of $h$ satisfying 
$m_e \in \Sigma$ and
\[\dim \Sigma=2(\dim G- \dim H).\]
\end{corollary}

\noindent\textbf{Proof} The manifold $\Sigma$ is the submanifold of  
the surface
$S$, obtained by setting the parameter $\alpha\in\fh$ equal to  
zero; in other 
words
\begin{equation}
\label{defining expression}
\Sigma=\{[g,\eta,v(\eta,0)]\in Y\mid g\in G,\eta\in\fm^*\}.
\end{equation}
The submanifold $\Sigma$ is a smooth manifold, since~(\ref{dimension
formula for relative}) implies that 
it has constant rank $2(\dim G- \dim H)$; that is, the map
\[\begin{array}{cccccc}
\mathcal{T}:&G\times\fm^*&\longrightarrow&G\times\fm^*\times
V_{m_e}&\longrightarrow&G\times_H(\fm^*\times V_{m_e})\\
 &(g,\eta)&\longmapsto&(g,\eta,v(\eta,0))&\longmapsto&[g,\eta,v(\eta,0)]
\end{array}\]
with image $\Sigma$ is a local constant rank map
around $(e,0)\in G\times\fm^*$
with rank equal to $2(\dim G- \dim H)$, which implies that the surface 
$\Sigma$ is locally a manifold through the relative equilibrium  
$m_e$, of 
dimension $2(\dim G- \dim H)$. (See, for instance, \cite[Theorem  
3.5.18]{mta}.) 

The symplectic nature of $\Sigma$ can be verified in a  
straightforward manner.
Indeed, we will check that if $i:\Sigma\hookrightarrow Y$ is the natural
inclusion then the pair $(\Sigma,\omega_\Sigma)$, with $\omega_\Sigma=
i^*\omega_Y$, is a symplectic submanifold of $(Y,\omega_Y)$.
Let $\pi:G\times\fm^*\times V_{m_e}\rightarrow G\times_H
(\fm^*\times V_{m_e})$ be the canonical projection. Note
that every  vector in $T_{[g,\eta,v(\eta,0)]}\Sigma$ can be
written as $T_{(g,\eta,v(\eta,0))}\pi(T_e L_g\cdot\zeta,\delta
\eta,D_{\fm^*}v(\eta,0)\cdot\delta)$, for some $\zeta\in
\fg$ and $\delta\eta\in\fm^*$. The two--form $\omega_\Sigma$
is clearly closed. In order to prove that it is nondegenerate, let's
suppose that the vector $T_{(g,\eta,v(\eta,0))}\pi(T_e L_g\cdot\zeta,
\delta\eta,D_{\fm^*}v(\eta,0)\cdot\delta \eta')$ is such that
\begin{eqnarray}
0 &=& \omega_\Sigma([g,\eta,v(\eta,0)])
(T_{(g,\eta,v(\eta,0))}\pi(T_e L_g\cdot\zeta,\delta\eta, \nonumber \\
&& \qquad \quad D_{\fm^*}v(\eta,0)\cdot\delta\eta),
T_{(g,\eta,v(\eta,0))}\pi(T_e L_g
\cdot\zeta',\delta\eta',D_{\fm^*}v(\eta,0)\cdot\delta\eta')) 
\label{big_expression}
\end{eqnarray}
for every $\zeta'\in\fg$ and $\delta\eta'\in\fm^*$.
We will show that this implies that $T_{(g,\eta,v(\eta,0))}\pi(T_e L_g
\cdot\zeta,\delta\eta,D_{\fm^*}v(\eta,0)\cdot\delta\eta)=0$ necessarily.
Using $\omega_\Sigma=i^*\omega_Y$ and the explicit expression of the
symplectic form $\omega_Y$ associated to the MGS normal form
(see the previously quoted original papers, as well as~\cite{rosou,  
re, thesis}), we
can write (\ref{big_expression}) in the form
\begin{eqnarray*}
0 &=& \langle\delta\eta' +
D\mathbf{J}_{V_{m_e}}(v(\eta,0))\cdot(D_{\fm^*}v(\eta,0)\cdot\delta\eta'),
\zeta\rangle-\langle\delta\eta + 
D\mathbf{J}_{V_{m_e}}(v(\eta,0))\cdot(D_{\fm^*}v(\eta,0)\cdot\delta\eta),  
\zeta'\rangle\\
&& \qquad + {}
\omega_{V_{m_e}}(D_{\fm^*}v(\eta,0)\cdot\delta\eta,
D_{\fm^*}v(\eta,0)\cdot\delta\eta')
\end{eqnarray*}
for any $\delta\eta'\in\fm^*$ and $\zeta'\in\fg$. If we fix
$\delta\eta'=0$ and let $\zeta'$ be arbitrary,
we obtain 
\[
\delta\eta +
D\mathbf{J}_{V_{m_e}}(v(\eta,0))\cdot(D_{\fm^*}v(\eta,0)\cdot\delta\eta)=0.
\]
Since $\delta\eta\in\fm^*$,
$D\mathbf{J}_{V_{m_e}}(v(\eta,0))\cdot(D_{\fm^*}v(\eta,0)\cdot\delta\eta)\in\fh^*$,
and $\fm^*\cap\fh^*=\{0\}$, we have 
\begin{equation}
\label{equal to 0}
\delta\eta=D\mathbf{J}_{V_{m_e}}(v(\eta,0))\cdot  
(D_{\fm^*}v(\eta,0)\cdot\delta\eta)=0.
\end{equation}
If we now fix $\zeta'=0$ and let $\delta\eta'$ be arbitrary,
we obtain $\zeta\in\fh$, which, together
with~(\ref{equal to 0}), guarantees that $T_{(g,\eta,v(\eta,0))}
\pi(T_e  
L_g\cdot\zeta,\delta\eta,D_{\fm^*}v(\eta,0)\cdot\delta\eta)=0$, as
required.   $\ \ \ \blacksquare$

\medskip

In the remainder of this section we will show that  the
persistence phenomena described by Theorem~\ref{abelian
bifurcation theorem} and Corollary~\ref{symplectic
manifold of equilibria} preserve stability. More
specifically, we will show that if the relative
equilibrium $m_e$ is stable, then the entire local
symplectic manifold $\Sigma$ given by
Corollary~\ref{symplectic manifold of equilibria}
consists of stable relative equilibria. First, we
recall the definition of nonlinear stability of a
relative equilibrium:

\begin{definition}
Let $(M,\omega,h,G,\mathbf{J}:M\rightarrow \fg^*)$ be a
Hamiltonian system with symmetry and let $G'$ be a subgroup of
$G$. A relative equilibrium $m_e\in M$ is called
$G'${\bfi --stable}, or {\bfi stable modulo} $G'$,
if for any $G'$--invariant open neighborhood $V$ of the
orbit $G'\cdot m_e$, there is an open neighborhood $U\subseteq V$ of
$m_e$, such that if $F_t$ is the flow of the Hamiltonian
vector field $X_h$ and $u\in U$, then $F_t (u)\in V$
for all $t\geq 0$.
\end{definition}

Before recalling the stability result that we will use in our discussion
we introduce the following notation. Suppose that we fix a
splitting of $\fg$ as in~(\ref{decompositions}). If
$\xi=\xi_1+\xi_2$, with $\xi_1\in\fg_{m_e}$   and $\xi_2\in\fm$,
is a generator of the relative equilibrium $m_e$, then the unique  
element
$\xi_2\in\fm$ is called the {\bfi orthogonal generator} of $m_e$ with
respect to the splitting~(\ref{decompositions}).

We now state the following theorem whose proof can be found  
in~\cite{singreleq}
or in~\cite{re}.

\begin{theorem}
\label{bd}
Let $(M,\,\{\cdot,\cdot\},\, h)$ be a Poisson system with a
symmetry given by the Lie group $G$ acting
properly on $M$ in a globally Hamiltonian fashion, with
associated equivariant momentum map $\mathbf{J} :M\rightarrow
\fg^*$. Assume that the Hamiltonian $h\in C^\infty (M)$ is
$G$--invariant. Let $m_e\in M$ be a relative equilibrium such that
$\mathbf{J}(m_e)=\mu\in\fg^*$, $\fg^*$ admits an ${\rm
Ad}^*_{G_\mu}$--invariant
inner product, $H:=G_{m_e}$, and $\xi\in
{\rm Lie}(N_{G_\mu}(H))$ is its orthogonal generator, relative to a
given ${\rm Ad}_H$--invariant splitting. If the quadratic form
\[
D^2(h-\mathbf{J}^\xi)(m_e)|_{W\times W}
\]
is definite for some (and hence for any) subspace $W$ such that
\[
\label{definitionofwre}
\ker D\mathbf{J}(m_e)=
 W\oplus T_{m_e}(G_\mu\cdot
m_e),
\]
then $m_e$ is a $G_\mu$--stable relative equilibrium.
If $\dim W=0$, then $m_e$ is
always a $G_\mu$--stable relative equilibrium. The quadratic
form $D^2(h-\mathbf{J}^\xi)(m_e)|_{W\times W}$, will be
called the {\bfi stability form}
of the relative equilibrium $m_e$.
\end{theorem}

A relative equilibrium satisfying the hypotheses
of Theorem~\ref{bd} is said to be {\bfi formally stable}. 
Note that in the
Abelian case all the adjoint invariance requirements in the statement of
the previous theorem are trivially satisfied. We now state our stability
persistence result.

\begin{proposition}
Under the conditions of Corollary~\ref{symplectic manifold of equilibria},
suppose that the relative equilibrium $m_e$ is  formally (and  
consequently
nonlinearly) stable; that is, it has an orthogonal generator $\xi\in\fm$
with respect to the splitting~(\ref{decompositions}) such that the  
quadratic form
\[
D^2(h-\mathbf{J}^\xi)(m_e)|_{W\times W}
\]
is definite for some (and hence for any) subspace $W$ such that
\[
\ker D\mathbf{J}(m_e)=
 W\oplus T_{m_e}(G\cdot
m_e),
\]
then the symplectic manifold $\Sigma$ of relative equilibria  
passing through
$m_e$
can be chosen (by taking, if necessary, a sufficiently
small neighborhood of $m_e$ in  the submanifold $\Sigma$ of  
Corollary~\ref{symplectic manifold of equilibria}) to consist  
exclusively of nonlinearly  stable relative
equilibria. 
\end{proposition}

\noindent\textbf{Proof} Recall that the symplectic manifold  $\Sigma$
consists of points of the form $[g,\eta,v(\eta,0)]\in Y$, with  
$\eta\in\fm^*$
sufficiently close to $0$, which are relative equilibria with  
generator $\xi+\beta
(\eta,0)$. Since $\xi\in\fm$ is by hypothesis an orthogonal  
generator with
respect to the splitting~(\ref{decompositions}), and the function  
$\beta$ maps
into $\fm$, the generator $\xi+\beta(\eta,0)$ is also an orthogonal
generator for the relative equilibrium $[g,\eta,v(\eta,0)]\in Y$.
Hence, in order to prove the Proposition it suffices to show that
the quadratic form
\[D^2(h-\mathbf{J}^{\xi+\beta(\eta,0)})([g,\eta,v(\eta,0)])|_{W_{[g,\eta,v(\eta,0)]}\times
 W_{[g,\eta,v(\eta,0)]}}\]
is definite for some subspace $W_{[g,\eta,v(\eta,0)]}$ such that $\ker
D\mathbf{J}([g,\eta,v(\eta,0)])=W_{[g,\eta,v(\eta,0)]}\oplus  
T_{[g,\eta,v(\eta,0)]}(G\cdot
[g,\eta,v(\eta,0)])$. Using the expression of the momentum map in the
MGS--coordinates
described in Proposition~\ref{marle tube}, it is easy to verify that
\begin{eqnarray*}
\lefteqn{\ker D\mathbf{J}([g,\eta,v(\eta,0)])} \\
&=&T_{[g,\eta,v(\eta,0)]}(G\cdot
[g,\eta,v(\eta,0)])\oplus  
T_{[e,\eta,v(\eta,0)]}\Phi_g(T_{(\eta,v(\eta,0))}
\Psi(\{0\}\times \ker D\mathbf{J}_{V_{m_e}}(v(\eta,0)))),
\end{eqnarray*}
where $\Phi_g$ denotes the $G$--action in MGS coordinates (see
Proposition~\ref{marle tube}) and $\Psi$ is the slice mapping
introduced in~(\ref{marle slice mapping}). This identity singles
out the space 
\[
T_{[e,\eta,v(\eta,0)]}\Phi_g(T_{(\eta,v(\eta,0))}
\Psi(\{0\}\times \ker D\mathbf{J}_{V_{m_e}}(v(\eta,0))))
\]
as a
choice for $W_{[g,\eta,v(\eta,0)]}$. We are now in position to
study the definiteness of the stability form of the relative
equilibrium $[g,\eta,v(\eta,0)]$, using as $W_{[g,\eta,v(\eta,0)]}$
the space just mentioned. Indeed,
\begin{align}
D^2(h&-\mathbf{J}^{\xi+\beta(\eta,0)})
([g,\eta,v(\eta,0)])|_{W_{[g,\eta,v(\eta,0)]}\times
W_{[g,\eta,v(\eta,0)]}}\notag\\
 &=D^2(h-\mathbf{J}^{\xi+\beta(\eta,0)})([g,\eta,v(\eta,0)])
 |_{(T_{[e,\eta,v(\eta,0)]}\Phi_g(T_{(\eta,v(\eta,0))}\Psi(\{0\}\times
 \ker D\mathbf{J}_{V_{m_e}}(v(\eta,0)))))\times ({\rm same})}\notag\\
 &=D^2((h-\mathbf{J}^{\xi+\beta(\eta,0)})\circ\Phi_g)([e,\eta,
 v(\eta,0)])|_{(T_{(\eta,v(\eta,0))}\Psi(\{0\}\times\ker  
D\mathbf{J}_{V_{m_e}}(v(\eta,0))))\times ({\rm same})}\notag\\
 &=D^2(\cH -\cJ^{\xi+\beta(\eta,0)})
 (\eta,v(\eta,0))|_{(\{0\}\times\ker  
D\mathbf{J}_{V_{m_e}}(v(\eta,0)))\times (\{0\}\times\ker  
D\mathbf{J}_{V_{m_e}}(v(\eta,0)))}.\label{almost there yes!}
\end{align}
The formal stability of $m_e$ implies that the quadratic form
\[
D^2\left(\cH - \cJ^{\xi}\right)(0,0)|_{(\{0\}\times
V_{m_e})\times(\{0\}\times V_{m_e})}
\]
is definite, therefore, since definiteness is an open condition, for any
$\eta\in\fm^*$ close enough to 0,
\[
D^2\left(\cH - \cJ^{\xi  + \beta(\eta,0)}\right)
(\eta,v(\eta,0))|_{(\{0\}\times V_{m_e})\times(\{0\}\times V_{m_e})}
\]
is also definite.
Since $\ker D\mathbf{J}_{V_{m_e}}(v(\eta,0))$ in
expression~(\ref{almost there yes!}) is a subset of $V_{m_e}$,
the definiteness of the stability form of relative equilibrium
$[g,\eta,v(\eta,0)]$
is guaranteed for small enough $\eta\in\fm^*$, as required.$\ \ \  
\blacksquare$

\section{Bifurcation of relative equilibria with maximal isotropy}
\label{abelian g mu and equilibria}

Consider the symmetric Hamiltonian
system $(M,\omega,h,G,\mathbf{J}:M\rightarrow\fg^{*})$,  where the  
Lie group $G$ acts properly 
on the manifold $M$.
Let $m_e \in M$ be a relative equilibrium with 
momentum $\mu = \bJ(m_e)$.  In contrast to the previous section,  
we will assume here that the 
relative equilibrium $m_e$ is degenerate, that is, there is a  
generator $\xi\in\g$ and  a 
nontrivial vector subspace $V_0\subset T_{m_e}M$ for which 
\begin{equation}
\label{degenerate hypothesis}
\ker \, D^2 (h - \bJ^\xi)(m_e) = \fg_{\mu} \cdot m_e \oplus V_0.
\end{equation}
This hypothesis implies that in writing the reduced critical point  
equations the Lyapunov--Schmidt reduction will be nontrivial and
there will be the possibility of genuine bifurcation.

In this section we will focus on the study of the bifurcation
equation~\textbf{(B1)}; that is, we will assume that 
the rigid residual equation is satisfied and therefore the
relative equilibria near $m_e$ correspond to the zeroes
of~\textbf{(B1)}. 

In the framework of general dynamical systems, the bifurcation of
relative equilibria with isotropy group $K$, out of a degenerate 
(i.e. nonhyperbolic) isolated equilibrium, is {\bfi generic} 
\footnote{Loosely speaking, a property of a system is generic if it  
is true
unless additional constraints are added to the system (see~\cite{g2}).}
if $K$ is maximal and satisfies an additional property, e.g. has an 
odd-dimensional fixed--point subspace in the space $V_0$ on which the
bifurcation equation is defined, or has an even dimensional fixed--point
subspace together with a non-trivial $S^1$ action. 
The famous Equivariant Branching Lemma (see, e.g., 
\cite{g2}), belongs to the former case, while the latter
appears in a work of Melbourne (see~\cite{melbourne 94, chossat  
maximal}).
We shall see that both results have a counterpart in the symmetric  
Hamiltonian 
case, although being Hamiltonian is a nongeneric property from the  
general dynamical systems point of view. 
When searching for relative equilibria, the generator ($\alpha\in\fg$) or 
momentum ($\eta\in\fg^*$) serves as a bifurcation parameter, in addition 
to any physical control parameters present in system. Due to the ``rigidity''
of these geometric ``parameters'', care must be taken when adapting the 
bifurcation theorems to relative equilibria of Hamiltonian systems.

As a final preliminary remark, we point out the fact that our  
theorems will 
be stated for bifurcation from a general relative equilibrium, not  
just from
one pure (isolated) equilibrium. In the latter case, the gradient  
character of
the bifurcation equation (see Remark~\ref{rem_gradient}) simplifies the 
arguments (see Remark~\ref{rem_simplifications}).

\subsection{A Hamiltonian Equivariant Branching Lemma} 

In the situation described above, let $m_e\in M$ be a relative  
equilibrium satisfying the 
degeneracy hypothesis~(\ref{degenerate hypothesis}). As we saw in 
Proposition~\ref{equivariance features}, the bifurcation  
equation~\textbf{(B1)} 
can be constructed so as to be $\sym$--equivariant, which implies  
that for any subgroup 
$K\subset\sym$, $B$ can be restricted to the $K$--fixed point  
subspaces in 
its domain and range; 
hence we can find solutions of $B$ by finding the solutions of 
\[
B^K:=B|_{(\fm^*)^K\times V_0^K\times\fg_{m_e}^K}:
(\fm^*)^K\times V_0^K\times\fg_{m_e}^K\longrightarrow V_2^K.
\]
Assume now that $K\subset\sym$ is a maximal isotropy subgroup of the 
$\sym$--action on $V_0$ and, moreover, that $\dim (V_0^K)=1$.
Under this hypothesis we will look for pairs 
$(\eta,v_0)\in (\fm^*)^K\times V_0^K$ satisfying
\begin{equation}
\label{alpha equal zero}
B^K(\eta,v_0,0)=0.
\end{equation}
Note that $\dim (V_0^K)=1$ implies that (see for instance~\cite{Bre}) 
\[
L:=N_{G_\xi\cap G_{m_e}}(K)/K\simeq\left\{
\begin{array}{l}
\{Id\}\\
	\\
\mathbb{Z}_2.
\end{array}
\right.
\]
Recall that $L$ acts naturally on $(\fm\sus)^K$ and on $V_0^K$, and  
that $B^K$ is $L$--equivariant. Depending on the character of the  
$L$--action, the first  terms in the Taylor expansion of~(\ref{alpha  
equal zero}) can be written as
\[
B^K(\eta,v_0,0)=\left\{
\begin{array}{lll}
\kappa\cdot \eta+ v_0^2 c+\cdots=0&\quad\text{if}\quad& L \simeq\{Id\}\\
	\\
v_0(\kappa\cdot \eta+v_0^2 c+\cdots)=0&\quad\text{if}\quad& L  
\simeq\mathbb{Z}_2,
\end{array}
\right.
\]
for some vector $\kappa\in(\fm\sus)^K$ and some constant $c$ that are 
generically nonzero.
In both instances these expressions generically allow us to solve for $v_0$ 
in terms of the other variables via the Implicit Function Theorem,  
giving us 
saddle--node type branches if $L\simeq\{id\}$
and a pitchfork bifurcation if $L\simeq\mathbb{Z}_2$ 
(see~\cite{g2} for arguments of this sort). More explicitly, we  
have proved 
the following result.

\begin{theorem}[Equivariant Branching Lemma]
\label{Equivariant Branching Lemma}
Let $m_e\in M$ be a relative equilibrium of the Hamiltonian system  
$(M,\omega,h,G,\mathbf{J}:M\rightarrow\fg^{*})$,  where the Lie  
group $G$ acts properly on the manifold $M$. Suppose that there is a  
generator $\xi\in\g$ and  a nontrivial vector subspace $V_0\subset  
T_{m_e}M$ for which 
\begin{equation*}
\ker \, D^2 (h - \bJ^\xi)(m_e) = \fg_{\mu} \cdot m_e \oplus V_0.
\end{equation*}
Then, generically, for any subgroup $K\subset G_\xi\cap G_{m_e}$  
for which 
$\dim (V_0^K)=1$ and the rigid residual equation is satisfied on 
$(\fm^*)^K\times V_0^K\times\{0\}$, a branch of relative equilibria  
with 
isotropy subgroup $K$ bifurcates from $m_e$. If 
$N_{G_\xi\cap G_{m_e}}(K)/K\simeq\{Id\}$, the bifurcation is a  
saddle--node;
if $N_{G_\xi\cap G_{m_e}}(K)/K\simeq\mathbb{Z}_2$, it is a pitchfork.
\end{theorem}
We will illustrate this result with an example in the following section.

\subsection{Bifurcation with maximal isotropy of complex type}

In what follows we will use a strategy similar to the one introduced by 
Melbourne~\cite{melbourne 94} in the study of general equivariant  
dynamical 
systems, to drop the hypothesis on the dimension of $V_0^K$ in the  
Equivariant 
Branching Lemma. Our setup will be the same as in  
Theorem~\ref{Equivariant 
Branching Lemma} but in this case we will be looking at maximal complex 
isotropy subgroups $K$ of the $\sym$--action on $V_0$, that is, maximal 
isotropy subgroups $K$ for which
\begin{equation}
\label{maximality}
L:=N_{\sym}(K)/K\simeq\left\{
\begin{array}{l}
S^1\\
	\\
S^1\times\mathbb{Z}_2.
\end{array}
\right.
\end{equation}
Notice that in such cases $V_0^K$ has even dimension.

As in the previous section, we will use the equivariance properties  
of the 
bifurcation equation in order to restrict the search for its  
solutions to the 
$K$--fixed space $(\fm^*)^K\times V_0^K\times\fg_{m_e}^K$.  
Moreover, we will 
consider only solutions of the form $(0,v_0,\alpha)\in 
(\fm^*)^K\times V_0^K\times\pe$, where $\pe$ is some
$\mbox{Ad}_{N_{\sym}(K)}$--invariant complement to $\mathfrak{k}$ in 
${\rm Lie}(N_{\sym}(K))$. Note that (\ref{maximality}) implies that
$\pe\simeq \mathfrak{l} \simeq \mathbb{R}$. 

We now show that the adjoint action of $N_{\sym}(K)$ on $\pe$ is  
trivial.
The canonical projection $\pi:N_{\sym}(K)\rightarrow L$ 
is a group homomorphism; hence the commutivity of $L$ implies that
\[
\pi \lp g h g^{-1} \rp = \pi(g) \pi(h) \pi(g)^{-1} = \pi(h)
\]
for any $g$, $h \in N_{\sym}(K)$. In particular, 
\[
T_e\pi\cdot(\mbox{Ad}_g \alpha)
= \ddto \pi \lp g \, \exp(t \, \alpha) \, g^{-1} \rp 
= \ddto \pi \exp(t\alpha) 
= T_{e}\pi\cdot\alpha
\]
for any $g \in N_{\sym}(K)$ and $\alpha \in {\rm  
Lie}(N_{\sym}(K))$, which
implies that $\mbox{Ad}_g - \mbox{id}$ maps ${\rm Lie}(N_{\sym}(K))$
into $\ker(T_{e}\pi)=\mathfrak{k}$. Since $\pe \cap \mathfrak{k}=\{0\}$ 
and $\pe$ is $\mbox{Ad}_{N_{\sym}(K)}$--invariant, it follows that 
$\lp \mbox{Ad}_g - \mbox{id} \rp|_{\pe} = 0$ for all $g \in  
N_{\sym}(K)$,
i.e. that the adjoint action on $\pe$ is trivial. 

\begin{theorem}
\label{maximal theorem}
Let $m_e\in M$ be a relative equilibrium of the Hamiltonian system 
$(M,\omega,h,G,\mathbf{J}:M\rightarrow\fg^{*})$,  where the Lie  
group $G$ acts 
properly on the manifold $M$. Suppose that there is a generator  
$\xi\in\g$ and  
a nontrivial vector subspace $V_0\subset T_{m_e}M$ for which 
\begin{equation*}
\ker \, D^2 (h - \bJ^\xi)(m_e) = \fg_{\mu} \cdot m_e \oplus V_0.
\end{equation*}
Suppose that the fixed point set $V_0^\sym=\{0\}$. Then for each maximal complex isotropy subgroup $K$ of the  
$\sym$--action 
on $V_0$ such that
\[
[{\rm Lie}(N_{\sym}(K)),\fg_{m_e}^K]=0
\]
and each $\mbox{Ad}_{N_{\sym}(K)}$--invariant complement $\pe$ to 
$\mathfrak{k}$ in ${\rm Lie}(N_{\sym}(K))$ such that the rigid residual 
equation $\rho(0, v_0, \alpha) = 0$ is satisfied for all $v_0 \in  
V_0^K$ 
and $\alpha \in \pe$,
there are generically at least $\frac{1}{2} \dim V_0^K$
(respectively $\frac{1}{4} \dim V_0^K$) branches of 
relative equilibria bifurcating from $m_e$ if 
$N_{\sym}(K)/K\simeq S^1$ (respectively $S^1\times\mathbb{Z}_2$).
\end{theorem}

\noindent\textbf{Proof} Let $B: \mathcal{U}_3 \subset 
\fm\sus\times V_0\times\fg_{m_e}\rightarrow V_2$ 
be the bifurcation equation corresponding to the reduced critical point 
equations constructed around $m_e$ using the MGS--slice mapping  
introduced 
in~(\ref{marle slice mapping}). The equivariance of this slice mapping 
guarantees that $B$ is $\sym$--equivariant; hence any solutions of 
\[
B^K:=B|_{(\fm^*)^K\times V_0^K\times\fg_{m_e}^K}:(\fm^*)^K\times V_0^K
\times\fg_{m_e}^K\longrightarrow V_2^K.
\]
are solutions of $B$.

As we stated above, we will restrict our search to
solutions in the set $\{0\}\times V_0^K\times\pe$, where
$\pe$ is some $\mbox{Ad}_{N_{\sym}(K)}$--invariant complement to 
$\mathfrak{k}$. Identify $V_0$ and $V_2$ using an invariant inner  
product and
define $\tilde B^K:V_0^K\times\pe\rightarrow V_0^K$ through the  
relations
\begin{equation}
\label{thing to deal with}
\langle \tilde B^K(v_0,\alpha), u \rangle 
:= B^K(0,v_0,\alpha) \cdot u
= D_{V_{m_e}} \lp \cH - \cJ^\eta \rp (0, v_0 + v_1 (0, v_0, \alpha)) 
\cdot u|_{\eta = \Xi(0,v_0,\alpha)}
\end{equation}
for any $v_0$, $u \in V_0^K$ and $\alpha \in \pe$. The equivariance 
properties of $B$ and the triviality of the action on $\pe$ imply that 
$\tilde B^K$ satisfies the following equivariance condition:
\begin{equation}
\label{case 1}
\tilde B^K(g\cdot v_0,\alpha)=g\cdot \tilde B^K(v_0,\alpha)
\quad\text{for all}\quad g\in N_{\sym}(K).
\end{equation}
Note that as a corollary to this property we have 
\begin{equation}
\label{0 at 0}
\tilde B^K(0,\alpha)=0\qquad\text{for all}\ \alpha,
\end{equation}
since for all $g\in N_{\sym}(K)$, $g\cdot \tilde B^K(0,\alpha)
=\tilde B^K(0,\alpha)$ and, consequently, the isotropy subgroup of 
$\tilde B^K(0,\alpha)$ contains $N_{\sym}(K)$ and hence it strictly contains 
$K$. By the maximality of $K$ as an isotropy subgroup, the isotropy 
subgroup of $\tilde B^K(0,\alpha)$ is necessarily $\sym$. However, given that 
by hypothesis $V_0^\sym=\{0\}$, we have $\tilde B^K(0,\alpha)=0$, as claimed.

We find the solution branches by first finding an open ball  
$B_r(0)$ about
the origin in $V_0^K$ and a function $\alpha: B_r(0)  
\rightarrow\pe$ satisfying
\[
\langle \tilde B^K(v_0,\alpha(v_0)), v_0 \rangle =0,
\]
then using $\tilde B^K$ and $\alpha$ to define a family of vector fields
on the unit sphere in $V_0^K$. Standard topological arguments show that
these vector fields have the requisite number of equilibria, which 
correspond to solutions of the original equations.

As the  first step in finding the function $\alpha$ we compute the  
Taylor expansion of $\tilde B^K$. As a result of the Lyapunov-Schmidt reduction and of~(\ref{0 at 0}), we can write
\[
\tilde B^K(v_0,\alpha)=L(\alpha)v_0 + g(v_0,\alpha),
\]
where $L(\alpha)$ is a linear operator such that $L(0)=0$, and  
$g(v_0,\alpha)$ is such that $g(0,\alpha)=0$, $D_{v_0}g(0,\alpha)=0$  
for all $\alpha$.
Moreover, a lengthy but straightforward computation shows that
\[
L(\alpha) = -\mathbb{P}D_{V_{m_e}V_{m_e}}\cJ^\alpha(0,0)
+ L_1(\alpha),
\]
where $L_1(0)=L_1'(0)=0$.
We now show that if we identify $V_0$ and $V_2$ by means of an  
invariant inner product, then there exists a constant  
$k\in\mathbb{N}^*$ such that 
\begin{equation}
\label{second derivative identity}
-\mathbb{P}D_{V_{m_e}V_{m_e}}\cJ^\alpha(0,0)|_{V_0^K}
= \alpha  \, k \, \mathbb{I}_{V_0^K},
\end{equation}
where $\mathbb{I}_{V_0^K}$ denotes the identity on $V_0^K$.
Indeed, note that 
\[
\cJ^\alpha(0,v)
=\langle\mathbf{J}\left(\Psi(0,v)\right),\alpha\rangle
=\langle\mu, \alpha \rangle + \mathbf{J}^\alpha_{V_{m_e}}(v)
\]
and hence
\begin{equation}
D_{V_{m_e}V_{m_e}}\cJ^\alpha(0,0)\cdot(v,w)
=D_{V_{m_e}V_{m_e}}\mathbf{J}_{V_{m_e}}^\alpha(0)(v,w)
=\omega_{V_{m_e}}(\alpha\cdot v,w)\label{second j}
\end{equation}
for any $v,\,w\in V_{m_e}$

We now restrict our attention to elements $v,\,w\in V_{m_e}^K$. 
Recall that since $V_{m_e}$ is symplectic, the vector subspace  
$V_{m_e}^K$ 
is symplectic with a canonical $L$ action; hence for any 
$\alpha\in\mathfrak{l}$ and $v\in V_{m_e}$ there is an infinitesimally 
symplectic transformation $A_\alpha$ such that $\alpha\cdot  
v=A_\alpha v$. 
The equivariant version of the Williamson normal form due to 
Melbourne and Dellnitz~\cite{md}, implies the existence of a basis  
in which 
$A_\alpha$ and $\omega_{V_{m_e}^K}$ admit simultaneous matrix  
representations
consisting of three diagonal blocks corresponding to the subspaces 
$E_\mathbb{R}$, $E_\mathbb{C}$, and $E_\mathbb{H}$ of $V_{m_e}^K$  
on which $L$ 
acts in a real, complex, and quaternionic fashion, respectively. 
Moreover, in this basis the restrictions of $A_\alpha$ and  
$\omega_{V_{m_e}^K}$ 
to $E_{\mathbb{C}}$ take the form:

\[
\omega_{V_{m_e}^K}|_{E_{\mathbb{C}}}= \pm i \, \mathbb{I} 
\qquad \mbox{and} \qquad
A_\alpha|_{E_{\mathbb{C}}}=
\pm i \, \alpha \, \mbox{diag}(k_1, \ldots, k_q).
\]
for some natural numbers $k_1,\ldots,\,k_q$. The signs in these two  
equalities are consistent, that is, they are either both positive  
or both negative (in all that follows we will focus only in the  
positive case). 
These expressions follow directly from the tables in~\cite{md} and the
absence of nilpotent parts in $A_\alpha$, which is dictated 
by the requirement that $A_\alpha$ be the zero matrix when $\alpha=0$. 
By hypothesis $K$ is a maximal isotropy subgroup of the 
$\sym$--action on $V_0$ for which $V_0^K\subset E_\mathbb{C}$.  
Moreover, 
since the $L$--action on $V_0^K\setminus\{0\}$ is free,
there exists $k\in\mathbb{N}^*$ such that
\[
A_\alpha|_{V_0^K}=  i \, k \, \alpha \, \mathbb{I}_{V_0^K}.
\]
Using this expression in~(\ref{second j}), we 
obtain~(\ref{second derivative identity}) and hence
\[
\tilde B^K(v_0,\alpha) = \alpha \, k \, v_0 + L_1(\alpha)v_0 +  
g(v_0,\alpha)
\]
where $L_1(\cdot)$ is of order higher than $|\alpha|$ and  
$g(\cdot,\alpha)$ is of order higher than $\|v_0\|$. It follows that  
in
the equation
\[
0 = \langle v_0,\tilde B^K(v_0,\alpha)\rangle
= \alpha k\|v_0\|^2 + \langle v_0,L_1(\alpha)v_0+g(v_0,\alpha)\rangle
\]
we can factor out $\|v_0\|^2$ and then apply the implicit function  
theorem to
obtain a unique function $\alpha: B_r(0) \to \pe$, for some $r > 0$,
near the solution $(0,0)$.

Using this function we can define a one parameter family of  
$L$--equivariant 
vector fields $X_\epsilon$ on $S^{2n - 1}$ by
\[
X_\epsilon(u) = \tilde B^K(\epsilon \, u, \alpha(\epsilon \, u)).
\]
The zeroes of these vector fields correspond to solutions of
the bifurcation equation. Since $L$ acts freely on $S^{2n-1}$,  
$X_\epsilon$ 
determines a smooth vector field $\tilde X_\epsilon$ on $S^{2n-1}/L$;
the Poincar\'e--Hopf theorem implies that $\tilde X_\epsilon$  
generically has 
at least  
\[
\chi\left(S^{2n-1}/L\right)
= \left \{ \begin{array}{ll}
\chi(\mathbb{CP}^{n-1})=n & \mbox{if $L\simeq S^1$} \\
\chi\left(\mathbb{CP}^{n-1}/\mathbb{Z}_2\right)=n/2 \qquad 
& \mbox{if $L\simeq S^1\times\mathbb{Z}_2$}
\end{array} \right .
\]
equilibria.

The following lemma proves that $X_\epsilon(u)$ is always  
orthogonal to the 
tangent space $\mathfrak{l} \cdot u$ of the $L$--orbit of $u$, 
i.e. $\langle X_\epsilon(u), \zeta_{S^{2n-1}}(u) \rangle = 0$ 
for any $u \in S^{2n-1}$ and $\zeta \in \mathfrak{l}$.
Hence the equilibria of $\tilde X_\epsilon$ correspond to orbits of  
equilibria 
of $X_\epsilon$, which in turn determine orbits of 
solutions of the bifurcation equation.

\begin{lemma}
\label{orthogonality}
If $[{\rm Lie}(N_{\sym}(K)),\g_{m_e}^K]=0$, then 
$\langle X_\epsilon(u), \zeta_{S^{2n-1}}(u) \rangle = 0$ 
for any $u \in S^{2n-1}$ and $\zeta \in \mathfrak{l}$.
\end{lemma}
\noindent\textbf{Proof} 
We first show that $\tilde B^K(v_0, \alpha)$ is orthogonal to 
$\mathfrak{l} \cdot v_0$ for any $v_0 \in V_0^K$ and $\alpha \in \pe$. 

Given $\alpha \in \pe$, define $\mathcal{H}_\alpha: V_0^K \to  
\mathbb{R}$ and
$\mathbf{j}_\alpha: V_0^K \to \fg^*$ by
\[
\mathcal{H}_\alpha(v_0) = \mathcal{H}(v_0 + v_1(0, v_0, \alpha))
\qquad \mbox{and} \qquad
\mathbf{j}_\alpha(v_0) = \mathbf{j}(v_0 + v_1(0, v_0, \alpha)).
\]
The equivariance of $v_1$ and triviality of the action on $\pe$ imply
that $\mathcal{H}_\alpha$ is $\sym$--invariant and  
$\mathbf{j}_\alpha$ is
$\sym$--equivariant. 

We can choose as a complement $V_1$ to $V_0$ in $V$ the space  
annihilated
by $V_2$. (If $V_2$ is identified with $V_0$ using an inner  
product, this 
choice for $V_1$ is the orthogonal complement to $V_0$ in $V$.) In  
this case,
\[
D_{V_{m_e}}\lp \mathcal{H}-\mathbf{j}^{\Xi(0, v_0, \alpha)} \rp
(0,v_0 + v_1(0, v_0, \alpha) )\cdot v_1 = 0 
\]
for any $v_0 \in V_0$, $v_1 \in V_1$, and $\alpha \in \fg_{m_e}$.  
Hence, 
given $v_0$, $u \in V_0^K$, $\alpha \in \pe$, and
$\zeta \in (\fg_{m_e} \cap \fg_\xi) \subset \fg_\mu$, if we set 
$\eta = \Xi(0,v_0,\alpha)$ and $v = v_0+v_1(v_0,\alpha)$, then
\begin{eqnarray*}
\langle \tilde B^K(v_0,\alpha), \zeta_{V_0}(v_0) \rangle
&=& D_{V_{m_e}}\lp \mathcal{H}-\mathbf{j}^\eta \rp(0,v)\cdot 
\zeta_{V_0}(v_0) \\
&=& D_{V_{m_e}}\lp \mathcal{H}-\mathbf{j}^\eta \rp(0,v)\cdot 
((\mbox{id} + D_{V_0} v_1(0, v_0, \alpha))\cdot \zeta_{V_0}(v_0)) \\
&=& D\lp \mathcal{H}_\alpha-\mathbf{j}_\alpha^\eta \rp(v_0)\cdot 
\zeta_{V_0}(v_0) \\
&=& \langle{\rm ad}\sus_{\zeta}\mathbf{j}_\alpha(v_0), \eta \rangle \\
&=& \langle {\rm ad}\sus_{\zeta}\mathbf{J}_{V_{m_e}}(v), \eta \rangle.
\end{eqnarray*}
In particular, if $\zeta \in {\rm Lie}(N_{\sym}(K))$ and
$[{\rm Lie}(N_{\sym}(K)),\g_{m_e}^K]=0$, then 
$\langle \tilde B^K(v_0,\alpha), \zeta_{V_0}(v_0) \rangle = 0$, 
since $\mathbf{J}_{V_{m_e}}(v)\in(\fg_{m_e}\sus)^K$.

To complete the proof, note that the linearity of the action  
implies that
\[
\langle X_\epsilon(u), \zeta_{S^{2n-1}}(u) \rangle
= \langle \tilde B^K(\epsilon \, u, \alpha (\epsilon \, u)), 
	\zeta_{S^{2n-1}}(u) \rangle
= \frac{1}{\epsilon} 
\langle \tilde B^K(\epsilon \, u, \alpha (\epsilon \, u)), 
	\zeta_{V_0}(\epsilon u) \rangle
= 0.
\]
\ \ \ \ $\blacksquare$

\begin{remark} 
\normalfont
Note that Theorem \ref{maximal theorem} provides a (generic) lower  
bound for 
the number of branches of critical points of bifurcating from  
$m_e$. In fact, 
if $m_e$ has nontrivial isotropy, then in many situations a sheet  
of critical 
points bifurcates from $m_e$, rather than a finite number of one  
dimensional 
branches. An example of this phenomenon is given in \S  
\ref{wave_section}.
A continuous curve of bifurcation points with nontrivial isotropy 
appears in many other symmetric Hamiltonian systems, including the  
Lagrange top 
and the Riemann ellipsoids. (See, for example, \cite{tops}, \cite{bld}, 
\cite{qpm}, and \cite{llt}.) In \cite{bld} it is shown that for  
Lagrangian
systems with $S^1$ symmetry this phenomenon occurs under conditions that
are generic within that class of systems. \ \ \ $\blacklozenge$

\begin{remark} \label{rem_simplifications}
\normalfont There are two cases in which the theorem that we just proved
applies in a particularly straightforward manner. First, suppose that
the relative equilibrium $m_e$ is such that its momentum value
$\mu=\mathbf{J}(m_e)$ has an Abelian isotropy subgroup $G_\mu$. In such
situation we have automatically that $[{\rm
Lie}(N_{\sym}(K)),\g_{m_e}^K]=0$ for any $K\subset\sym\subset G_\mu$ and
also, using the techniques introduced in Section~\ref{Slice map
refinements and the rigid residual equation} (see especially
Corollary~\ref{abelian matters}), condition~\textbf{(ii)} on the rigid
residual equation can be easily dealt with.

Another case of interest is when $m_e$ is actually a pure  
equilibrium whose
isotropy is the entire symmetry group $G$, i.e. the $G$--orbit of $m_e$ 
reduces to $m_e$ itself. Note that in that case
$\fm=\fp=\{0\}$ and therefore the rigid residual equation just does not
exist. Also, the condition $[{\rm Lie}(N_{\sym}(K)),\g_{m_e}^K]=0$ in
the statement of the theorem is not necessary in that case since 
the bifurcation equation is variational (see  
Remark~\ref{rem_gradient}) and
therefore the associated vector field is orthogonal to the $G$--orbits, 
and {\em a fortiori} to the $N_{\sym}(K))$--orbits in $V_0^K$. 
It is interesting to note that in this case, the Equivariant  
Branching Lemma
stated in Theorem~\ref{Equivariant Branching Lemma} is not applicable, 
because the parameter $\eta$ is now missing.\ \ \ $\blacklozenge$
\end{remark}

\section{An example from wave resonance in mechanical systems}
\label{wave_section}

When a  Hamiltonian mechanical system has two natural modes with
frequencies in the  approximate ratio $p:q$, $p$ and $q$ integers, one
says it has a $p:q$ internal resonance. When in addition the system
possesses  rotational as well as reflectional invariance, and the
natural modes are wavy, that is, they break the rotational  
invariance in the form
of $m$ waves for mode $p$, say, and $n$ waves for mode $q$, then the
reflectional invariance implies the coexistence of 4 complex amplitudes
associated to these modes, namely $z_1e^{mi\varphi}$ and
$z_2e^{-mi\varphi}$ for mode $p$, and $z_3e^{ni\varphi}$ and
$z_4e^{-ni\varphi}$ for mode $q$. The case with $(p,q)=(m,n)=(1,2)$ is
of particular interest, for example, in the analysis of certain water wave
problems in a cylindrical geometry, and has been studied with particular
attention in~\cite{ChoDias}. In this work, the method of investigation
was the projection of the vector field {\em in normal form} on the orbit
space of the group action which here is $G=O(2)\times S^1$, where the
$S^1$ term comes from the normal form assumption. The simple
form of the group action allowed, by this method, a fairly detailed
description of the possible dynamics near the trivial state. In  
particular, various kinds of relative equilibria were described in  
this way. 
This is therefore a good example to test our approach to the bifurcation
analysis of relative equilibria.

A straightforward change of coordinates in the group $SO(2)\times S^1$
allows us to define the toral part of the action of $G$ as
\begin{equation}
(\phi,\psi) \mapsto  
(z_1e^{i\phi},z_2e^{i\psi},z_3e^{2i\phi},z_4e^{2i\psi}),
~~(\phi,\psi)\in S^1\times S^1,
\end{equation}
where $(z_1,z_2,z_3,z_4)\in \mathbb{C}^4$ are the complex  
amplitudes of the  
system. 
The reflection acts by permutation of $z_1$ with $z_2$ and of $z_3$ with
$z_4$. Then, it was shown in~\cite{Bridges} that the general form of
a $G$-invariant, real smooth Hamiltonian $h$ is
$$
h = h(X_1,X_2,X_3,X_4,U_1,U_2)
$$
where 
\begin{eqnarray}
X_j &=& z_j\bar z_j \\
U_k &=& \frac{1}{2}(z_k^2\bar z_{k+2} + {\bar z}_k^2z_{k+2}),~~k=1,2.
\end{eqnarray}
and $h$ is invariant under the (simultaneous) permutation of $X_1$ with
$X_2$, $X_3$ with $X_4$ and $U_1$ with $U_2$. 

The Lie algebra $\g$ of $G$ acts on $\mathbb{C}^4$ by
\begin{equation}
(\xi_1,\xi_2)\in \mathbb{R}^2 \mapsto (i\xi_1z_1,i\xi_2z_2,2i\xi_1z_3,2i\xi_2z_4)
\end{equation}
and the momentum map $\J$ is defined by
\begin{equation}
\J(z_1,z_2,z_3,z_4)= \left( \begin{array}{c} |z_1|^2+2|z_3|^2 \\ 
|z_2|^2+2|z_4|^2 \end{array}\right),
\end{equation}
where we have identified $\g^*$ with $\mathbb{R}^2$.

We now write the  relative equilibrium equation, $D(h-\J^\xi)(m) =0$  in
complex coordinates. We set $\xi=(\xi_1,\xi_2)$ and
$$a_j=\frac{\partial h}{\partial X_j},~~
b_k=\frac{\partial h}{\partial U_k}. $$
Note that these coefficients are real and satisfy the relations
$$a_{2k}(X_1,X_2,X_3,X_4,U_1,U_2)= a_{2k-1}(X_2,X_1,X_4,X_3,U_2,U_1),~~ 
k=1,2,$$ 
and 
$$b_2(X_1,X_2,X_3,X_4,U_1,U_2)=b_1(X_2,X_1,X_4,X_3,U_2,U_1).$$ 
Thus 
\begin{eqnarray}
\lefteqn{D \lp h - J^\xi \rp(z)} \nonumber \\
&&= \lp (a_1 - \xi_1) z_1 + b_1\bar z_1z_3,
(a_2 - \xi_2) z_2 + b_2\bar z_2z_4,
(a_3 - 2 \xi_1) z_3 + \frac{b_1 z_1^2}{2},
(a_4 - 2 \xi_2) z_4 + \frac{b_2 z_2^2}{2} \rp.
\label{wave_eqn}
\end{eqnarray}

We can use the symmetries of the system (\ref{wave_eqn}) to easily  
identify 
a branch of relative equilibria; we will then use the results of the 
previous sections to find other branches of relative equilibria  
bifurcating 
from this branch. Note that the subgroup $H=\mathbb{Z}_2\times S^1$ 
of the toral group is an isotropy subgroup of $G$,
with fixed-point subspace $\mbox{Fix}(H)=\{(0,0,z_3,0)~|~z_3\in  
\mathbb{C}\}$.
This space is invariant under the map $D \lp h - J^\xi \rp$;  
specifically,
\[
D \lp h - J^\xi \rp(0, 0, z_3, 0) = (0, 0, (a_3 - 2\xi_1) z_3, 0).
\]
(Here we make the obvious identification
of $\mathbb{C}^4$ and $\lp \mathbb{C}^4 \rp^*$.) 
Thus every element of $\mbox{Fix}(H)$ is a relative equilibrium,  
each with
a one--parameter family of generators $(\hat \xi_1, \xi_2)$, where
\begin{equation} \label{xi1}
\hat\xi_1 = {\textstyle \frac{1}{2}} a_3(0,0,X_3,0,0,0)
\end{equation}
and $\xi_2$ is arbitrary.
Note that $(\hat\xi_1,0)$ is an orthogonal generator, while
$\fh = \lcb (0, \alpha) : \alpha \in \mathbb{R} \rcb$. 
The trajectory of each such relative equilibrium is 
$$
z(t) = (0, 0, Ce^{i\hat\xi_1 t+\varphi}, 0)
$$
and is parameterized by a positive number $C$ and phase $\varphi$. 
We call this family of relative equilibria $RE_I$
and analyze the bifurcation of new relative equilibria from this
family by applying our slice map decomposition at the points
$z_e=(0,0,C,0)$. Recall that, the isotropy subgroup of $z_e$ is 
$H=\mathbb{Z}_2\times S^1$. 

In constructing a slice mapping using Proposition~\ref{build_slice},
note that the linearity of the phase space $\mathbb{C}^4$ allows us  
to use 
the trivial chart map $\psi(u)=z_e+u$, where $u\in \mathbb{C}^4$. 
The linearization of the momentum map at $z_e$ is
\[
D \J(z_e) \cdot (\delta z_1, \delta z_2, \delta z_3, \delta z_4) 
= \lp \begin{array}{c} 4 C \, \mbox{Re}(\delta z_3) \\ 0  
\end{array} \rp,
\]
with $\mbox{ker} \, D \J(z_e) = 
\lcb (z_1, z_2, i y, z_4) : z_j \in \mathbb{C}, y \in \mathbb{R} \rcb$.
We set $\fm$ equal to the orthogonal complement to $\fh = \fg_{z_e}$ in 
$\fg_\mu = \fg$ and set $V$ equal to the orthogonal complement
to $\fg \cdot z_e = \{(0,0,2 i \xi C,0)~|~\xi\in \bR\}$ in 
$\mbox{ker} \, D \J(z_e)$, so that
\[
\fm \approx \fm^* = \{(\eta ,0)~|~\eta \in\bR\}
\qquad \mbox{and} \qquad
V=\{(z_1,z_2,0,z_4) | z_j \in \bC \}. 
\]
Finally, we set $W=\{(0,0,\eta,0)~|~\eta\in\bR\}$.
These choices yield the slice map
\[
\Psi(\eta, v) := \lp 0, 0, n(\eta), 0 \rp + v
= \lp z_1, z_2, n(\eta), z_4 \rp,
\quad \mbox{where} \quad n(\eta) := C + \frac \eta {4 C}.
\]
The pullback of the energy--momentum function by the slice map $\Psi$ is
\[
(\cH - \cJ^\xi)(\eta, v) = h \lp X_1, X_2, n^2, X_4, n \, Y_1, U_2 \rp
	- \xi_1 \lp X_1 + 2 n^2 \rp - \xi_2 (X_2 + 2 \, X_4),
\]
where $Y_1 := \mbox{re}(z_1^2)$ and $n = n(\eta)$.

The analysis of the relative equilibria is simplified by the  
commutativity of
$\fg$,
which implies that $\g_\mu=\g$ and the two ``rigid''  
equilibrium 
conditions ({\bf RE1}) and ({\bf RE2}) are trivially satisfied.  
Hence the
first nontrivial step in the algorithm is Step 2: The map  
$\omega_1$ is simply
$\omega_1(\eta, \beta, \alpha) = (\hat \xi_1 + \beta, \xi_2 +  
\alpha)$. Hence 
solving 
\[
0 = D_{\fm^*}(\cH - \cJ^{\omega_1})(\eta, v) 
= \frac {b_1 Y_1 - 4 \beta \, n} {4 C} 
\]
for $\beta$ yields $\beta(\eta, v) 
:= \frac {b_1 Y_1}{4 n(\eta)}$. 
Thus 
\[
\Xi(\eta, v, \alpha) = (\hat \xi_1 + \beta(\eta, v), \xi_2 + \alpha)
= \lp \frac {a_3} 2 + \frac {b_1 Y_1}{4 n(\eta)}, \xi_2 + \alpha \rp. 
\]
and 
\begin{eqnarray}
\lefteqn{D_V(\cH - \cJ^{\Xi(\eta, v, \alpha)})(\eta, v) } \nonumber \\
&& = \lp 2 \lp a_1 - {\textstyle \frac {a_3} 2} - \beta \rp z_1 
	+ b_1 n \bar z_1, 
	2 (a_2 - (\xi_2 + \alpha)) z_2 + b_2 \bar z_2 z_4, 
	2 (a_4 - 2 (\xi_2 + \alpha)) z_4 + {\textstyle \frac {b_2  
z_2^2} 2} \rp.
\label{bif_eq_wave}
\end{eqnarray}

The bifurcation of relative equilibria from $(RE_I)$ depends on the 
invertibility of the linearization of the relative equilibrium equation 
in $V$ at the point $(0,0)$. The second variation 
$D_{VV} \lp \cH - \cJ^\Xi \rp (0, 0)$ has eigenvalues and eigenspaces
\[
\begin{array}{lcllcl}
\lambda_1^+ &=& a_1 - {\textstyle \frac {a_3} 2} + C b_1
\qquad & V_1^+ &=& \{(x, 0, 0, 0) \ | \ x \in \bR \} \\
\lambda_1^- &=& a_1 - {\textstyle \frac {a_3} 2} - C b_1
\qquad & V_1^- &=& \{(i \, y, 0, 0, 0) \ | \ y \in \bR \} \\
\lambda_2 &=& a_2 - \xi_2 \quad \mbox{(double)} \qquad &
V_2 &=& \{(0, z_2, 0, 0) \ | \ z_2 \in \bC \} \\
\lambda_4 &=& a_4 - 2 \xi_2 \quad \mbox{(double)} \qquad &
V_4 &=& \{(0, 0, 0, z_4) \ | \ z_4 \in \bC \}.
\end{array}
\]
Note that the isotypic decomposition of $V$ with respect to the  
action of $H$
guarantees the decomposition of $D_{VV} \lp \cH - \cJ^\Xi \rp (0,  
0)$ 
into three $2\times 2$ blocks associated to $V_1^+ \oplus V_1^-$,  
$V_2$, and
$V_4$, since the action of $S^1$ separates the $z_1$ component from
the $z_2$ and $z_4$ ones, while the action of $Z_2$ separates  
further the
$z_2$ component from the $z_4$ one. 

There are two kinds of bifurcation points: 
\begin{enumerate}
\item {\em Bifurcation at $\lambda_1^+$ or $\lambda_1^-=0$}. 
Because these are simple eigenvalues, the Lyapunov-Schmidt  
procedure yields
one-dimensional bifurcation equations. The conditions of the Hamiltonian
Equivariant Branching Lemma are met, hence we can conclude the  
existence of a 
bifurcated branch of relative equilibria parameterized by $\eta\in  
\bR$ at 
each of these points. 
Note that ${\mathbb Z}_2$ acts as $- \mbox{Id}$ on the eigenvectors
associated with these eigenvalues. Thus it follows that the  
bifurcation is
of {\em pitchfork} type. The isotropy group of these solutions still
contains $S^1$. Therefore these relative equilibria fill 1-tori, i.e.
are still periodic solutions for the Hamiltonian vector field. 

Note that in this case, $C$ can be taken as the bifurcation parameter.
This is however the same thing as  taking $\eta$, since $W$ is  
defined as 
the subspace $\{(0,0,C+\eta,0)\}$ in $\mathbb{C}^4$. 

\item {\em Bifurcation at $\lambda_2=0$ or $\lambda_4=0$}. 
In both of these cases, the eigenvalue is double and therefore the  
space 
$V_0$ determined by the Lyapunov--Schmidt procedure is two dimensional.
Note that $SO(2)$ acts nontrivially on $\ker(A_2)$ and  
$\ker(A_4)$. Therefore 
the isotropy subgroup is maximal of complex type in both cases. Applying
Theorem~\ref{maximal theorem} yields at least one 
branch of circles of relative equilibria in each case. In fact, there 
is a two--parameter family (modulo symmetry) of relative equilibria
containing ($RE_I$).
These solutions live on 2-tori and are quasi-periodic whenever the
ratio of the two components of the generator is irrational.
What distinguishes these two families, aside from the fact that they
bifurcate at different values of $\xi_2$, is their symmetry: the  
isotropy
of the solutions bifurcating in the $z_4$ direction is $\mathbb{Z}_2$,
while it reduces to the trivial group for those bifurcating in the  
$z_2$ 
direction.
\end{enumerate}
Note that while the bifurcations associated to $\lambda_1^\pm = 0$  
generically 
occur only at isolated values of $C$, the bifurcations associated  
to $\lambda_2 = 0$ 
and $\lambda_4 = 0$ occur for any value of $C$ satisfying the  
nondegeneracy condition
$a_2(0, 0, C^2, 0) \neq a_4(0, 0, C^2, 0)$, since the second  
component $\xi_2$
of the generator at $z_e$ can always be chosen to equal $a_2(z_e)$ or 
$a_4(z_e)$. 

We now proceed with the actual solution of the bifurcation equation.
We first consider the bifurcation at $\lambda_1^+=0$. Generically
the remaining eigenvalues are nonzero at this point; we shall consider
only this case. Since $V_1^+$ is invariant under $D_V(\cH - \cJ^\Xi)$,
uniqueness of $v_1$ implies that $v_1 \equiv 0$ and the bifurcation
equation {\bf (B1)} is simply $D_V(\cH - \cJ^\Xi)|_{V_1^+} = 0$, i.e.
\[
0 = D_V(\cH - \cJ^{\Xi(\eta, (x_1, 0, 0), \alpha)})(\eta, (x_1, 0, 0)) 
= (f_1(\eta, x_1^2) x_1, 0, 0),
\]
where
\[f_1(\eta, s) := 2 a_1(s, 0, n^2, 0, n \, s, 0) 
	- a_3(s, 0, n^2, 0, n \, s, 0)
	- {\textstyle \frac s n}b_1(s, 0, n^2, 0, n \, s, 0),
	\qquad n = n(\eta). 
\]
Unless we are in the highly degenerate case $D_\eta f_1(0, 0) = 
D_s f_1(0, 0) = 0$, we can use the Implicit Function Theorem to solve 
for one variable in terms of the other. If, for example, we solve for
$\eta$ as a function of $s$, we obtain a unique function 
$\eta: (-\epsilon, \epsilon) \to \bR$ for some $\epsilon > 0$ satisfying
$f_1(\eta(s), s) = 0$, and hence 
$D_V(\cH - \cJ^{\Xi(\eta(x_1^2), (x_1, 0, 0), \alpha)})
(\eta(x_1^2), (x_1, 0, 0)) = 0$ for all $x_1^2 \in [0, \epsilon)$.  
Implicit 
differentiation of $f_2(\eta(s), s) = 0$ yields 
$\eta(s) = \frac s {4C} + o(s^2)$. 
Note that the group $\{0\}\times S^1$ is an isotropy subgroup of  
$G$, with
fixed-point space ${z_2=z_4=0}$. The bifurcation  
under consideration takes place in this subspace. 
The case $\lambda_1^- = 0$ is entirely analogous.

We now consider the case $\lambda_1^\pm \neq \lambda_2= 0 \neq  
\lambda_4$.
Application of the Lyapunov--Schmidt procedure yields 
\[
v_1(\eta, X_2) = (0, 0, 0, z_4(\eta, X_2)),
\qquad \mbox{where} \qquad 
z_4(\eta, X_2) := \frac {b_2 z_2^2}{2(a_4 - 2(a_2 + \alpha))}.
\]
Substituting $v_1$ into $D_V(\cH - \cJ^{\Xi(\eta, v, \alpha)})(\eta, v)$
yields 
\[
B(\eta, z_2, \alpha) = 
D_V(\cH - \cJ^{\Xi(\eta, z_2 + v_1(\eta, X_2), \alpha)})
	(\eta, z_2 + v_1(\eta, X_2))
= (0, f_2(\eta, X_2, \alpha) z_2, 0, 0),
\]
where
\[
f_2(\eta, X_2, \alpha) := \frac {b_2^2 X_2}{2(a_4 - 2 (a_2 + \alpha))}
	+ \alpha;
\]
here $a_2$, $a_4$, and $b_2$ are all evaluated at $(0, X_2,  
n(\eta)^2, 0, 0)$.
Since $f_2(0, 0, 0) = 0$ and $D_\alpha f_2(0, 0, 0) = 1$, there exists a
neighborhood $\mathcal{W}$ of $(0, 0)$ in $\bR \times [0, \infty)$ 
and a function $\alpha: \mathcal{W} \to \bR$ such that
$f_2(\eta, X_2, \alpha(\eta, X_2)) = 0$ for all $(\eta, X_2) \in  
\mathcal{W}$.
Since $f_2$ depends on $z_2$ only through $X_2 = |z_2|^2$, each  
zero of $f_2$
determines a circle of critical points of $\cH - \cJ^\Xi$. 
The case $\lambda_4 = 0$ is entirely analogous.

Note that in the cases $\lambda_2 = 0$ and $\lambda_4 = 0$, varying  
the parameter $\eta$ 
simply shifts the real component of $z_3$, and hence is equivalent  
to shifting
the initial relative equilibrium $z_e = (0, 0, C, 0)$; thus, when  
computing the 
complete bifurcation diagram near the line $\{ (0, 0, C, 0): C \in  
\bR \}$, 
we find that generically two pitchforks of revolution, one  
corresponding to 
$\lambda_2 = 0$ and the other to $\lambda_4 = 0$, emerge from each  
point $(0, 0, C, 0)$. 
In addition, there may be isolated points at which conventional (one 
dimensional) pitchforks emerge, corresponding to $\lambda_1^\pm = 0$.
\ \ \ $\blacklozenge$
\end{remark}

\end{document}